\documentclass{article}
\usepackage{amsfonts,graphicx}

\newtheorem{dfn}{Definition}[subsection] 
\newtheorem{prp}[dfn]{Proposition} 
\newtheorem{lmm}[dfn]{Lemma} 
\makeatletter

\@addtoreset{equation}{section}
\makeatother
\def \N {{\mathbb N}}
\def \Z {{\mathbb Z}}
\def \R {{\mathbb R}}

\def \Q {{\mathbb Q}}
\newcommand{\cst}{{\rm cst\,}}
\newcommand{\one}{{\tt 1 \hspace{-0.8ex} \tt l}}
\newcommand{\qed}{\hspace*{2ex} \hfill $\Box$}

\newcommand{\grad}{{\rm grad}}
\renewcommand{\div}{{\rm div}}
\newcommand{\X}{{\cal X}}
\newcommand{\E}{{\cal E}}
\newcommand{\D}{{\cal D}}
\renewcommand{\L}{{\cal L}}
\renewcommand{\P}{{\mathbb P}}
\newcommand{\U}{{\cal U}}
\newcommand{\Var}{{\rm Var}}

\title{Condenser physics applied to Markov chains\\
{\large -- A brief introduction to potential theory --}}
\author{A. Gaudilli\`ere}
\date{Ouro Preto, August 2008}

\begin{document}

\maketitle

These notes are based on a collective discussion that started with
Alessandra Bianchi in December 2007 at the WIAS in Berlin and went
on during the next two months in Rome 2 and 3 with Tony Iovanella,
Francesco Manzo, Francesca Nardi, Koli Ndreca, Enzo Olivieri,
Betta and Benedetto Scoppola, Alessio Troiani and Massimiliano
Viale. Since all of us were working on arguments more or less
strongly related to metastability, we all felt the need to
understand in some or more depth the links between potential theory and
Markov processes on which are based the tools coming from the
former that Bovier, Eckoff, Gayrard and Klein introduced in the
study of metastability \cite{BEGK1} and were successfully applied
in \cite{BEGK2}, \cite{BM}, \cite{BEGK3}, \cite{BdHN}, \cite{BBI} among other
papers. Many mathematicians and physicists are certainly well
aware of these links but, at various degrees, it was not our case and we thought that a
good understanding of this connection was an essential support to
intuition and a precious guide in the use of these tools. As J.L.
Doob once put it: {\em ``To learn potential theory from
probability is like learning algebraic geometry without
geometry.''}~\cite{Do} This explains why before going to the specific application
to metastability I will essentially focus on potential theory and Markov
processes for themselves. 

Still because concerned by the leading ideas founding the
connections between potential theory and Markov processes I will
write complete proofs only in the simplest context that allows for
avoiding any technicality and makes more transparent these
connections: that of Markov chains on a finite state space $\X$
(and this does not exclude working in the regime where the
cardinality $|\X|$ goes to infinity). We will not however
restrict our analysis to this simpler setting. Potential theory
begun indeed in a quite different form during the last three
decades of the XVIII century with the works of Lagrange and
Laplace that described the gravitational field as deriving from a
{\em potential} $V$ solution of the Laplace equation
\begin{equation}
\Delta V \equiv 0
\end{equation}
and blossomed in the first half of the XIX century
as the cornerstone of electrostatic,
in particular with Gauss and Green's works.
This development was so important
that the study of harmonic functions,
that is of the solutions of the Laplace equation,
used to be
one of the main pillar of the accademic
formation
of any physicist or mathematician
at the end of the century.

The study of Markov processes started with Markov at the beginning
of the next century but, as far as I know, it was not before the
last years of the second world war that the links with potential
theory begun to be drawn by Kakutani~\cite{Ka}. It took then
almost 40 years to make this connection fully developed. In 1984
Doob published his treaty {\em Classical Potential Theory and its
Probabilistic Counterpart} \cite{Do} and in the same year Doyle
and Snell wrote their beautiful article \cite{DS} that embraced in
a same light the mathematics of random walks and the physics of
electrical networks. Since then people did not stop harvesting the
fruits of such a fertile union. See for example: \cite{Ly},
\cite{GKZ}, \cite{Be}, \cite{ABS}, \cite{BF}, \cite{BdHS}.

A final motivation for writing these notes is that we could not
find (although it probably exists) a single synthetic text that
linked together the electrostatic of original potential theory,
the physics of electric networks and the probabilistic meaning of
the objects they contemplate. But I want to point out a few
classics that I found particularly useful
to write these notes,
even though some are
not mainly or directly linked to the subject: together
with Doyle and Snell's article~\cite{DS} there were Norris' book
\cite{No}, Lyons and Peres' continuously updated online book
\cite{LP}, Karatzas and Shreve's book \cite{KS}, Sinclair's paper
\cite{Si}, and Lawler's book \cite{La}.
I want also to thank here
Pietro Caputo and Alessandra Faggionato
for the many discussions we had
that considerably enriched these notes.

\section{Laplace equation}
\subsection{Harmonic functions}

In this section and the next one $ \U$ will denote an open subset
of the $d$-dimensio\-nal vector space $\R^d$.
For $x$ in $\R^d$ we will write $x_1$, $x_2$,~\dots, $x_d$
for its coordinates and we will use the notation
$B_2(x,r)$ for the open Euclidean ball of centre $x$ and radius $r>0$. 

\begin{dfn}
  For $f\in {\cal C}^2( \U,\R)$
  we define the {\em Laplacian}
  of $f$ as 
  \begin{equation}
  \Delta f : x\in  \U \mapsto \Delta_x f := \sum_{k=1}^{d} \frac{\partial^2 f}{\partial x_k^2} (x) \in \R
  \end{equation}
  and we say that $f$ is {\em harmonic} on $ \U$
  when it satisfies the {\em Laplace equation on~$ \U$} 
  \begin{equation}
  \forall x \in  \U,\: \Delta_x f =0
  \end{equation}
\end{dfn}

\noindent
{\bf Examples:} With
\begin{equation}
r: (x_1,\dots,x_d) \in \R^d \mapsto \sqrt{x_1^2+\dots+x_d^2}
\end{equation}
a harmonic function on $\R^d\setminus\{0\}$
is $r$ if $d=1$, $\ln r$ if $d=2$ and $r^{2-d}$ if $d\geq 3$.

\medskip\par\noindent
The previous definition is a good one
in the measure in which it frames harmonic functions
inside differential calculus
with all the tools it provides.
But this definition looks to harmonic functions
from an essentially local point of view.
It has to be reformulated
to make transparent larger scale properties
of harmonic functions.

For $f\in {\cal C}^1( \U,\R)$ we will denote by $\nabla f$
or $\vec{\nabla f}$ or $\grad f$ its gradient
\begin{equation}
\grad f : x\in \U
\mapsto
\nabla_x f:= \left(\frac{\partial f}{\partial x_k}(x)\right)_{1\leq k\leq d} \in \R^d
\end{equation}
For $\phi \in {\cal C}^1( \U, \R^d)$ we will denote by $\div\phi$
or $\nabla\phi$ its divergence
\begin{equation}
\nabla\phi : x\in \U \mapsto \div_x \phi := \sum_{k=1}^d \frac{\partial \phi_k}{\partial x_k}(x) \in \R
\end{equation}
Hence, saying that $f$ is harmonic on $ \U$ is saying that
\begin{equation}
\div(\grad f) \equiv 0
\end{equation}
or, equivalently, that $\nabla f$ is a null divergence field. 
Now Stokes' lemma makes possible to switch from the local point of view 
to a larger scale one.

\medskip\par\noindent
{\bf Lemma [Stokes]}
\it
  Let $\vec{\phi} \in {\cal C}^1( \U,\R^d)$
  and ${\cal V}$ be an open subset of $ \U$
  with a compact closure $\overline{{\cal V}} \subset  \U$
  and such that $\partial {\cal V}$ is a (smooth) submanifold of $\R^d$.
  Then
  \begin{equation}
  \int_{\partial {\cal V}} \vec{\phi} . \vec{dS} 
  := \int_{\partial {\cal V}} \vec{\phi} . \vec{n} \, d\sigma
  =\int_{\cal V} \div \vec{\phi} \, d\lambda
  \end{equation}
  where $\vec{n}$ is the unitary vector that is orthogonal to $\partial {\cal V}$ and oriented
  from ${\cal V}$ towards  $ \U\setminus {\cal V}$,
  while $\sigma$ and $\lambda$ stand for the surface and volume Lebesgue measure respectively.
\rm

\medskip\par\noindent
Stokes' lemma is a straightforward  identity
in its discrete version (see Section~\ref{flwsandcrrnts})
and we just assume it in its continuous one.
It implies that for a harmonic function $f$ on $ \U$
and such a closed oriented smooth surface $\partial {\cal V} \subset  \U$
the flux of the vector field $\nabla f$ through $\partial {\cal V}$ is zero.  
And as a first consequence we get:

\begin{prp}[Mean-value property]
  If $f$ is harmonic on $ \U$
  then $f$ satisfies the {\em mean-value property} (m.v.p.),
  that is:
  \begin{equation}
  \forall r>0,\, \forall x\in  \U,\,
  \overline{B_2(x,r)} \subset  \U
  \Rightarrow
  f(x) =  \int_{\partial B_2(x,r)} f \frac{d\sigma}{|\partial B_2(x,r))|}
  \end{equation}
  where $|\partial B_2(x,r)|$ denotes the surface area of $\partial B_2(x,r)$,
  in such a way that the integral computes the mean value
  of $f$ on $\partial B_2(x,r)$.
\end{prp}

\noindent
{\bf Proof:} We pick $x$ in $ \U$ and define,
for all $r>0$  such that $\overline{B_2(x,r)}\subset \U$,
\begin{equation}
g(r) :=   \int_{\partial B_2(x,r)} f\, \frac{d\sigma}{|\partial B_2(x,r)|}\\
= \int_{\partial B_2(0,1)} f(x+ru) d\sigma_1(u)
\label{dfng}
\end{equation} 
with $\sigma_1$ 
the uniform probability measure on
$\partial B_2(0,1)$.
By continuity of $f$,
\begin{equation}
\lim_{r\rightarrow 0} g(r) = f(x)
\end{equation}
Then, we just need to show that $g$ is constant,
i.e., has a null derivative.
Fix $r$ small enough to have $g(r)$ that is well defined.
For all small enough real $h$
and $u\in\partial B_2(0,1)$, the Taylor formula gives
\begin{equation}
f(x+(r+h)u) = 
  f(x+ru) + h \frac{\partial f}{\partial u}(x+ru)
  +\int_0^h (h -t) \frac{\partial^2 f}{\partial u^2}(x+(r+t)h) dt
\end{equation} 
so that,
integrating over $\partial B_2(0,1)$
and using $f\in {\cal C}^2( \U)$
to control the second derivative in the integral,
\begin{equation}
g(r+h)= g(r)
+ \frac{h}{|\partial B_2(x,r)|}
\int_{\partial B_2(x,r)} \vec{\nabla f}.\vec{dS} + o(h) 
\end{equation}
By Stokes' lemma the integral in this sum is zero
and we conclude $g'(r) = 0$.
\qed

\medskip\par
Actually the mean-value property characterizes harmonic functions
and gives additional information on their regularity:

\begin{prp}
  If $f\in {\cal C}( \U)$ has the m.v.p. then
  \begin{itemize}
  \item[i)]
  $f \in {\cal C}^\infty ( \U)$;
  \item[ii)] f is harmonic on $ \U$.
  \end{itemize}
\end{prp}

\noindent
{\bf Proof:} The proof of i) is based on a simple convolution argument
that can be find in \cite{KS} page 242.
Now, if $f$ is not harmonic, then we can find $x \in  \U$
and $r_0>0$ such that $\div \nabla f$ is strictly positive
(or strictly negative) on $B_2(x,r_0)\subset \U$.
Then the derivation of the previous proof
shows, with Stokes' lemma, that the function $g$ defined in~(\ref{dfng})
is strictly monotone in the neighbourhood of $0$.
And this contradicts the m.v.p. \qed

\medskip\par
The m.v.p. leads also to the following

\begin{prp}[Maximum principle]
  If $f$ is harmonic on $ \U$ then, for all compact sets
  $K\subset\overline{ \U}$ such that $f$ can be extended
  by continuity on $K$, $f|_K$ reaches its maximum (and its minimum)
  on $\partial K$.
\end{prp}

\noindent
{\bf Proof:} Set $M:= \max f(K)$. If $\partial K$ and $f|_K^{-1}\{M\}$
were disjoint sets then we could find $x$ in $\partial f|_K^{-1}\{M\}$
and $r>0$ such that $\overline{B_2(x,r)} \subset K\cap \U$.
The m.v.p would give
\begin{equation}
M=f(x)=\int_{\partial B_2(x,r)} f \frac{d\sigma}{|\partial B_2(x,r)|}
\end{equation}
while the continuity of $f$ in some $y\in(\partial B_2(x,r))\setminus f|_K^{-1}\{M\}$
would give   
\begin{equation}
\int_{\partial B_2(x,r)} f \frac{d\sigma}{|\partial B_2(x,r)|}<M
\end{equation}
what would be a contradiction. \qed

\medskip\par
The maximum principle is sometimes reported in the electrostatic
context as: 
{\em ``The potential has no local extremum
where there is no charge.''}
Indeed Gauss' law in Maxwell's equations reads
\begin{equation}
\div\vec{E}=\frac{\rho}{\epsilon_0}
\label{max}
\end{equation}
where $\rho$ stands for the charge density,
$\epsilon_0$ is the electric constant
and $\vec{E}$ is the electric field
that derives from a potential $V$, that is
\begin{equation}
\vec{E}=-\nabla V
\label{dfnV}
\end{equation}
A ``local extremum where there is no charge''
would be an isolated local extremum
somewhere in the interior $ \U$
of $\rho^{-1}\{0\}$, where the potential $V$ is harmonic
by~(\ref{max}) and~(\ref{dfnV}). 
And that would be in contradiction with the maximum principle.

The maximum principle opens the door to uniqueness
properties of the solution of the Dirichlet problem.

\begin{dfn}[Dirichlet problem]
Given $g\in {\cal C}(\partial  \U, \R)$ we say 
that $f$ is a solution of the Dirichlet problem
on $ \U$ with boundary condition $g$
if $f$ in ${\cal C}^2( \U)$ and ${\cal C}(\overline{ \U})$
satisfies the Laplace equation on $ \U$ and coincides
with $g$ on $\partial  \U$.
\end{dfn}

\noindent
{\bf Examples:} {\bf i)}
Consider a compact thermal conductor
$K$, fix to $g(x)$ the temperature in each point $x$
of $\partial K$ and assume that $g$ is continuous
on $\partial K$. If the temperature reaches
an equilibrium $f(x)$ in each point of the interior $ \U$
of the conductor, then $f$ will be solution of the Dirichlet problem
on $ \U$ with boundary condition $g$.
Let us assume the existence of such an equilibrium temperature.
The maximum principle gives us the uniqueness
of the equilibrium temperature. Indeed if $f_1$ and $f_2$
are both solutions, then $f_1-f_2$ is solution of the Dirichlet
problem with zero boundary condition. Then, on $ \U$, $f_1 -f_2$
cannot take values larger than the maximum value on the border,
that is 0, or smaller that the minimum value on the border,
0 once again: $f_1$ and $f_2$ coincide both on $ \U$ and~$\partial  \U$.
\smallskip\par\noindent
{\bf ii)} Consider a finite number of (disjoint) compact electric conductors
$A_1$,~\dots, $A_n$ 
and fix at values $V_1$,~\dots, $V_n$
on these conductors the difference of potential
with infinity. There cannot be any charge outside
the conductors, so that, by~(\ref{max}), (\ref{dfnV})
and taking the convention that the potential is $0$ at infinity,
a potential $V$ has to be solution of the Dirichlet
problem on the complementary of $\cup_n A_n$
with boundary condition $V_k$ on $\partial A_k$ for $k$ in
$\{1;\dots;n\}$ and with the additional condition
\begin{equation}
\lim_{x\rightarrow \infty} V(x) = 0
\end{equation}
Once again the maximum principle
gives the uniqueness of the potential $V$
under an existence hypothesis. 

\medskip\par
Proving the existence of a solution of a Dirichlet
problem turns out to be a rather difficult task
when one stay inside the framework of plain functional analysis.
It is time to turn to Markov processes.

\subsection{Brownian motion}

For simplicity we will now assume
that $ \U$ is a bounded open domain.
For extensions and generalizations to unbounded domains
of the results presented here we refer to~\cite{KS}
section~4.2. 
We denote by $P_x$ the law of a $d$-dimensional
Brownian motion $W$
starting from $x\in \R^d$
and by $\tau_A$ the hitting time of any set~$A$:
\begin{equation}
\tau_A=\inf\left\{t\geq 0 :\: W(t)\in A\right\}
\end{equation} 
Katunani's idea \cite{Ka} was to present the candidate
\begin{equation}
h: x\in \overline \U\mapsto E_x\left[g\left(W(\tau_{\partial \U})\right)\right]
\label{dfnh}
\end{equation}
to solution of the Dirichlet problem on $ \U$
with boundary condition $g$.
Since $\overline \U$ is a compact set
we have, for all $x$ in $\overline \U$, 
\begin{equation}
P_x(\tau_{\partial \U}<+\infty)=1
\end{equation}
so that $h$ is well defined.
We clearly have $h|_{\partial  \U} = g$
and $h$ has the m.v.p. Indeed, for any $x\in \U$
and $r>0$  such that $\overline{B_2(x,r)}\subset \U$,
we have, by 
the strong Markov property
at time $\tau_{\partial B_2(x,r)}$
and  using 
radial symmetry: 
\begin{eqnarray}
h(x)
&=& E_x\left[
  g\left(
    W(\tau_{\partial \U})
  \right)
\right]\\
&=& E_x\left[
  E_x\left[
    g\left(
      W(\tau_{\partial \U})
    \right)\Big| W(\tau_{\partial B_2(x,r)})
  \right]
\right]\\
&=& \int_{y\in\partial B_2(x,r)} 
E_y\left[
  g\left(
    W(\tau_{\partial \U})
  \right)
\right]
dP_x(W(\tau_{\partial B_2(x,r)})=y)\\
&=& \int_{\partial B_2(x,r)} 
h(y) \frac{d\sigma(y)}{|\partial B_2(x,r)|}
\end{eqnarray}
As a consequence $h$ is harmonic on $ \U$
and the only point to check
to get a solution of the Dirichlet problem
is the continuity of $h$ on $\overline \U$.
This question is intimately linked to the notion
of {\em regularity}.

\begin{dfn}
  For any set $A$ we define
  \begin{equation}
  \tau^+_A:=\inf\left\{t>0 :\: W(t) \in A\right\}
  \end{equation}
  and we say that {\em $ \U$ has a regular border} when
  \begin{equation}
  \forall a\in\partial \U,\; P_a\left(\tau^+_{ \U^c} = 0\right) =1
  \end{equation}
\end{dfn}

\begin{prp}
  A bounded open domain $\U$ has a regular border
  if and only if,
  for all $g$ in ${\cal C}(\partial \U)$
  the function $h$ defined in~(\ref{dfnh})
  is continuous on $\overline \U$,i.e.,
  is solution of the associated Dirichlet problem.
\end{prp}

We refer to \cite{KS} for the proof (of a stronger result).
We give now some examples.

In dimension two and for $b>a>0$
consider the Dirichlet problem on 
\begin{equation}
 \U=B_2(0,b)\setminus \overline{B_2(0,a)}
\end{equation}
with boundary conditions 
$1$ on $\partial B_2(0,a)$
and
$0$ on $\partial B_2(0,b)$.

On the one hand, since $\ln r$
is harmonic, we have the solution 
\begin{equation}
f=\frac{\ln b - \ln r}{\ln b - \ln a}
\end{equation}
By the maximum principle
$f$ is the unique solution
($\overline \U$ is a compact set).

On the other hand $ \U$
has a regular border.
Indeed a Brownian motion that starts
from $x$ in $\partial  \U$
crosses $\partial  \U$ infinitely
many times during any time interval
$[0;t]$ with $t>0$.
As a consequence the function $h$
defined in~(\ref{dfnh}) is solution 
of the problem and coincides with $f$.
That reads
\begin{equation}
\forall x \in \R^2,\;
a\leq r(x) \leq b \Rightarrow
P_x\left(\tau_{\partial B_2(0,a)}<\tau_{\partial B_2(0,b)}\right)
= \frac{\ln b - \ln r(x)}{\ln b - \ln a}
\label{lnpot}
\end{equation}

Consider now the Dirichlet problem on
the punctured ball
\begin{equation}
 \U = B_2(0,b)\setminus \{0\}
\end{equation}
This example will have some relevance later
when dealing with metastability in large volume
(Section \ref{metal}).
Sending $a$ to 0 in~(\ref{lnpot})
we get
\begin{equation}
\forall x \in \overline{B_2(0,b)}\setminus\{0\},\;
P_x\left(\tau_{\{0\}}<\tau_{\partial B_2(0,b)}\right)
= 0
\label{polar}
\end{equation}
This implies that the function $h$
defined in~(\ref{dfnh}) is equal to $\one_{\{0\}}$.
It is not continuous and 0 is not regular,
as it could directly be seen from~(\ref{polar}).
Such examples of domain with a non regular border
can be built with a connected border in dimension
$d\geq 3$.

\medskip\par
In the last example our candidate to solution for our problem
on the punctured ball lost the election.
But could have we find another solution?
The answer is no: using the uniqueness
of an eventual solution $f$ and the radial
symmetry of the problem we can show
that $f$ would be a simple function of $r$,
then solving the Laplace equation in polar coordinates
we would get 
\begin{equation}
f = \alpha + \beta\ln r
\end{equation}
with $\alpha$ and $\beta$ constants,
the continuity in 0 would imply $\beta = 0$
and the boundary conditions
could not be conciliated.
More generally:

\begin{prp}
  If the Dirichlet problem on
  a bounded open domain $ \U$
  with boundary condition $g$ has a solution $f$,
  then $f$ coincides with $h$ defined in~(\ref{dfnh}).
\end{prp}

\noindent
{\bf Proof:}
For all positive integer $n$ we define
\begin{equation}
 \U_n:=\left\{
  x\in  \U :\:
  d_2(x,\partial \U)> \frac{1}{n}
\right\}
\end{equation}
Note that $ \U_n$ has a regular border
since for all $a$ in $\partial \U_n$
$a$ is on the border of a ball (of radius 1/n)
contained in $ \U_n^c$.
We also define
\begin{equation}
h_n(x): x\in \overline \U \mapsto E_x\left[f(W(\tau_{ \U_n^c}))\right]
\end{equation}
The functions $f$ and $h_n$ coincide on $ \U_n^c$,
and they coincide on $\U_n$ too:
since $\U_n$ has a regular border, both are solution
of the Dirichlet problem on $\U_n$ with 
boundary condition $f|_{\partial \U_n}$,
since $\overline{\U_n}$ is a compact set,
this solution is unique by the maximum principle.
Now $f$ is bounded as continuous function
on the compact set $\overline \U$,
and by dominated convergence
we get, for all $x$ in $\overline \U$,
\begin{equation}
f(x)= \lim_{n\rightarrow +\infty} h_n(x) = h(x)
\end{equation}
\qed

The probabilistic approach to potential theory
does not only solve some of the problems
regarding the existence of a potential.
It also laid the ground
to receive deep insight from potential theory
into Markov processes theory.
For example formula~(\ref{lnpot})
gives the recurrence of the two-dimensional
Brownian motion: send $b$ to infinity to get
\begin{equation}
\forall x\not\in B_2(0,a),\; P_x(\tau_{\partial B_2(0,a)}<+\infty) = 1
\end{equation}
The same potential study in dimension $d\geq 3$
gives 
\begin{equation}
\forall x \in \R^d,\;
a\leq r(x) \leq b \Rightarrow
P_x\left(\tau_{\partial B_2(0,a)}<\tau_{\partial B_2(0,b)}\right)
= \frac{r(x)^{2-d}-b^{2-d}}{a^{2-d}-b^{2-d}}
\label{powpot}
\end{equation}
and the transience of the Brownian motion:
\begin{equation}
\forall x\not\in \overline{B_2(0,a)},\; P_x(\tau_{\partial B_2(0,a)}<+\infty)
= \left(\frac{r}{a}\right)^{2-d} < 1
\end{equation}

We close this section with a last illustration
of the evocative power of Kakutani's solution.
Consider a single electric compact conductor $K$
at potential 1 in $\R^3$.
By the so-called ``point-effect''
the electric field will be stronger
in the neighbourhood of the points $a$
of the convex parts of $K$
with strong curvature.
Indeed 
\begin{equation}
\vec{E}=-\nabla V
\end{equation}
and Kakutani's solution
for the potential $V$ outside $K$
gives
\begin{equation}
V(x) = P_x\left(
  \tau_K < +\infty
\right)
\end{equation}
(take the increasing limit of the solution
of the Dirichlet problem on $B_2(0,R)\setminus K$
with potential 0 on $\partial B_2(0,R)$
when $R$ goes to infinity.)
Loosely speaking 
(we will give a stronger justification
of the point-effect in Section~(\ref{cnc}))
the escape probability
can decrease much faster in the neighbourhood
of such points $a$.
And this why lightning rods are rods:
a stronger field 
in the neighbourhood of the rod makes easier
to reach the disruptive field
there before than somewhere else. 

\subsection{Discrete Laplacian and simple random walks}

Going from $\R^d$ to $\Z^d$
we loose the differential tool:
derivatives have to be replaced
by their discrete version.
Denoting by $\E$ the symmetric subset of $\Z^d\times\Z^d$
made of the nearest neighbour sites
\begin{equation}
\E := \left\{e=(e_-,e_+) \in (\Z^d)^2 :\: d_1(e_-, e_+)=1\right\}
\end{equation}
the gradient of a real valued function $f$ on $\Z^d$
becomes an antisymmetric real valued function on $\E$:
\begin{equation}  
\nabla f : e\in\E\mapsto\nabla_e f:= f(e_+)-f(e_-) \in \R
\end{equation}
and the divergence of such an antisymmetric function $\phi$
on $\E$ turns to be a real valued function on $\Z^d$:
\begin{equation}
\div \phi : x \in \Z^d \mapsto \div_x\phi
:=\sum_{\stackrel{\scriptstyle e \in \E}{\scriptstyle e_-=x}} \phi(e)\in \R
\end{equation}
The discrete Laplacian of $f$ is then defined by
\begin{equation}
\Delta f: x\in\Z^d\mapsto\Delta_x f:=\div_x(\nabla f) 
=\sum_{\stackrel{\scriptstyle y \in \Z^d}{\scriptstyle d_1(x,y)=1}} (f(y)-f(x))
\end{equation}
Note that this is coherent with
second order Taylor developments:
for $f$ in ${\cal C}^2(\R^d)$ and a unitary vector $u$
\begin{equation}
\frac{\partial^2 f}{\partial u^2} (x)=
\lim_{\stackrel{\scriptstyle h\rightarrow 0}{\scriptstyle h\in \R}}
\frac{f(x+hu)+f(x-hu)-2f(x)}{h^2}
\end{equation}
For $\U \subset \Z^d$ the border of $\U$ is
\begin{equation}
\partial \U := \left\{ e = (e_-,e_+) \in \E :\: e_-\in \U, e_+\not\in \U \right\}
\end{equation}
and its external border is
\begin{equation}
\partial_+ \U := \left\{ e_+ \in \Z^d :\: \exists e=(e_-,e_+) \in \partial \U\right\}
\end{equation}
A function $f$ defined on $\U\cup\partial_+ \U$ is {\em harmonic} on $\U$ if
\begin{equation}
\forall x\in \U,\; \Delta_x f = 0
\label{dfnhrmnc}
\end{equation}
Observe that~(\ref{dfnhrmnc}) expresses a local mean-value property:
it is equivalent to 
\begin{equation} 
\forall x\in \U,\; \frac{1}{2d}\Delta_x f 
=\left(
  \frac{1}{2d}\sum_{\scriptstyle d_1(x,y)=1} f(y)
\right)-f(x)
=0
\end{equation}
The set of harmonic functions on $\Z^d$
is the kernel of the generator
of the (continuous time) simple random walk,~$\frac{1}{2d}\Delta$, just like
the set of harmonic functions on $\R^d$ was the kernel
of the generator of Brownian motion,~$\frac{1}{2}\Delta$.

Like in the continuous case the mean-value property
of harmonic functions gives a maximum principle
(for which the notion of compactness
is replaced by that of finiteness),
and this maximum principle can be used
to show uniqueness properties for the solutions of Dirichlet problems.
Given $\U\subset\Z^d$ and $g$ a real valued function
on its external border,
we say that $f$ is a solution of the Dirichlet
problem on $\U$ with boundary condition $g$
if $f$ is harmonic on $\U$ and $f$ coincides with $g$
on $\partial_+ \U$.
Just like we used Brownian motion to prove the existence
of a solution for some Dirichlet problem,
we can do the same with simple random walks
on $\Z^d$. For example:

\begin{prp}
  For any finite subset $\U$ of $\Z^d$
  and any real valued function $g$ on $\partial_+ \U$,
  there is a unique solution of the Dirichlet problem on $\U$ with boundary
  condition $g$. This solution is the function $h$ defined by
  \begin{equation}
  \forall x \in \U\cup\partial_+ \U,\;
  h(x):= E_x\Big[
    g\left(
      \zeta(\tau_{\U^c})
    \right)
  \Big]
  \end{equation}
  where $E_x$ stands for the expectation
  under the law of a simple random walk~$\zeta$
  that start from $x$.
\end{prp}

\noindent
{\bf Proof:}
By the Markov property,
$h$ satisfies the local mean-value property on~$\U$.
Since $h$ and $g$ coincide on $\partial_+ \U$,
$h$ is solution of the Dirichlet problem.
It is the only one
by application of the maximum principle.
\qed

\section{Electrical networks}

\subsection{Random walks and generators}

An {\em electrical network} is a connected undirected 
weighted graph with positive weights,
with no more than one edge between
any pair of vertices and 
with finite total weight on each vertex.
More formally it is a pair $(\X, c)$
with $\X$ a countable set and $c$ a real valued non-negative
symmetric function on $\X\times\X$
such that
\begin{eqnarray}
&\forall x \in \X,\; \mu(x):= \sum_{y\in\X} c(x,y) < +\infty&
\end{eqnarray}
and such that, for all distinct $x$ and $y$ in $\X$,
there exist $x=z_1$, $z_2$,~\dots, $z_n=y$ in $\X$
with
\begin{equation}
\forall k \in \{1;\dots;n-1\},\; c(z_k,z_{k+1})>0
\end{equation}  
We call {\em nodes} the elements of $\X$,
we say that two nodes $x$ and $y$ are connected
when $c(x,y) >0$ and we call {\em edges}
the elements of $\E$, defined as the
set of ordered pairs of connected nodes:
\begin{equation}
\E:=\left\{
  (x,y)\in\X\times\X :\: c(x,y)>0
\right\}
\end{equation}
Of course the edges in $\E$ do have a direction,
but it is better to keep in mind
the image of an undirected graph
for which each pair 
of connected nodes can have
two representatives in the symmetric subset $\E$
of $\X\times\X$.
The {\em conductance} between two nodes $x$
and $y$ is $c(x,y)$ and the {\em resistance}
between $x$ and $y$ is
\begin{equation} 
r(x,y) := \frac{1}{c(x,y)} \in ]0;+\infty]
\end{equation}
Note that 0 and $+\infty$ are possible values for $c(x,y)$ 
and $r(x,y)$: when $(x,y)$ is not in $\E$.

We call {\em potential} any real valued function
on $\X$.
If we impose a potential $g(x)$ on each
node $x$ outside a subset $\U$ of $\X$,
an {\em equilibrium potential} $V$ associated
with the constraint
\begin{equation}
\forall x \in \U^c,\; V(x)= g(x)
\end{equation}
has to satisfy Ohm's and Kirchoff's laws.

\noindent
{\bf Ohm's law:}
The {\em current} $i$ associated with $V$ is
\begin{equation}
i: (x,y)\in\E\mapsto i(x,y)=\frac{V(x) - V(y)}{r(x,y)}
\end{equation}

\noindent
{\bf Kirchoff's law:}
For all $x$ in $\U$
\begin{equation}
\sum_{\stackrel{\scriptstyle y\in\X}{\scriptstyle (x,y)\in\E}}
i(x,y)
=0
\end{equation}
In other words
\begin{equation}
\forall x \in \U,\; -\L_x V = 0
\label{eq}
\end{equation}
with, for any potential $f$,
\begin{equation}
\L f : x\in\X\mapsto\L_x f:= \sum_{y\in\X}\frac{c(x,y)}{\mu(x)}(f(y)-f(x))
\label{dfnl}
\end{equation}
The operator $\L$ is the generator of $\xi$,
discrete time random walk on the network 
with transition probabilities 
\begin{equation}
p(x,y)=\frac{c(x,y)}{\mu(x)} \qquad x,y\in\X
\label{dfnp}
\end{equation}
(we call generator of a Markov chain
the generator of the associated continuous time Markov process
that updates its position at each ring of a Poissonian clock of intensity 1
according to the transition probabilities of the Markov chain)
and~(\ref{eq}) expresses once again
a local mean-value property
(that is also a martingale property
for the process $f(\xi)$ stopped
in $\U^c$).
As a consequence one can deal 
with the question of existence
and uniqueness of an equilibrium potential
associated with $\U$ and $g$
by using the maximum principle
that follows from the m.v.p.
and using Kakutani's solution.
For example if $\X$ is finite there exists
a unique equilibrium potential 
\begin{equation}
V(x)=E_x\left[g(\xi(\tau_{\U^c}))\right], \quad x\in\X
\end{equation}

\medskip\par\noindent
{\bf Remarks:}
{\bf i)} The Markov chain $\xi$ we associated
with $(\X,c)$ is ergodic and reversible with respect
to the measure $\mu$.
Conversely, {\em any}
reversible ergodic Markov chain $\xi$ on $\X$
is the random walk associated 
with some electrical network on $\X$.
If $\mu$ is a reversible measure
and $p(.,.)$ gives the transition probabilities
of $\xi$ we just define $c$ through~(\ref{dfnp})
to build a network for which the transition
probabilities of the associated random walk
are given by $p(.,.)$.
\smallskip\par\noindent
{\bf ii)} With each network $(\X,c)$ are associated
a unique random walk $\xi$ and, for each $\U\subset\X$,
a unique set ${\cal H}_\U$
of harmonic functions on $\U$, that is of solutions
of~(\ref{eq})-(\ref{dfnl}).
But an ergodic reversible random walk $\xi$ 
is associated with more than one
network since its associated reversible measure
is defined up to a multiplicative constant only.
These different networks correspond to different
choices of the conductance unity. Of course
when $\xi$ is associated with a {\em finite}
reversible measure there is a canonical choice
for the conductance unity: that for which $\mu$
is a probability.

A given family of sets of harmonic functions ${\cal H}_\U$
is associated with many more networks.
Indeed, from an electrical point of view the diagonal
values $c(x,x)$ of a network $(\X,c)$ are irrelevant
(note that self-loops are possible according
to our definitions).
Two electrical networks that differ only in these diagonal
values give rise to the same sets of harmonic potentials
${\cal H}_\U$,
but they are associated with quite different random walks
that do not have the same reversible measures.

\subsection{Flows and currents} \label{flwsandcrrnts}

For any $e = (x,y) \in \E$
and any potential $f$ 
we will use the notation
\begin{eqnarray}
e_- &=& x\\
e_+ &=& y\\
-e &=& (y,x)\\
\nabla_e f &=& f(y)-f(x) 
\end{eqnarray}
A {\em path} $\gamma$ is a finite or infinite sequence of edges
$e^1$, $e^2$,~\dots\ such that, for all $e^k$ and $e^{k+1}$ in $\gamma$
it is
\begin{equation}
e^k_+=e^{k+1}_-
\end{equation}
A {\em cycle} $\bar\gamma$ is a finite path $e^1$,~\dots, $e^n$
such that
\begin{equation}
e^n_+=e^1_-
\end{equation}  
We call {\em flow} any antisymmetric
real valued function on $\E$. 
The current $i=-c\nabla f$ associated, by Ohm's law,
with any potential $f$ is a flow.
But not all the flows derive from a potential.
It is easy to see that a flow $\phi$ derives from a potential
if and only if it satisfies the

\noindent
{\bf Second Kirchoff's law:}
\em
  For all cycle $\bar\gamma$
  \begin{equation}
  \sum_{e\in\bar\gamma} \phi(e) = 0
  \end{equation}
\rm

\noindent
In this case the associated potential
is uniquely defined up to an additive constant.

\medskip\par
The {\em divergence} of any flow $\phi$
is defined by
\begin{equation}
\div\phi: x\in\X\mapsto\div_x\phi:=\sum_{e_-=x}\phi(e)\in\R
\end{equation}
The border of any $\U\subset\X$ is
\begin{equation}
\partial \U := \left\{e\in\E :\: e_-\in \U, e_+\not\in \U\right\}
\end{equation}
and we  have

\medskip\par\noindent
{\bf Lemma [Stokes]:}
\em
  For any flow $\phi$ and any finite $K\subset\X$  
  \begin{equation}
  \sum_{e\in\partial K} \phi(e) = \sum_{x\in K} \div_x\phi
  \end{equation}
\rm

\medskip\par\noindent
{\bf Proof:}
\begin{eqnarray}
\sum_{x\in K} \div_x\phi
&=& \sum_{x\in K} \sum_{e_- = x} \phi(e)\\
&=& \sum_{\stackrel{\scriptstyle e\in\partial K}{\scriptstyle e_-\in K}} \phi(e)
+ \sum_{\stackrel{\scriptstyle e\not\in\partial K}{\scriptstyle e_-\in K}} \phi(e)
\end{eqnarray}
We call $S$ the last sum. It is also equal to
\begin{equation}
\sum_{\stackrel{\scriptstyle -e \not\in\partial K}{\scriptstyle (-e)_-\in K}} \phi(e)
= \sum_{\stackrel{\scriptstyle -e \not\in\partial K}{\scriptstyle (-e)_-\in K}} -\phi(-e)
\end{equation}
so that
\begin{equation}
S=-S=0
\end{equation}
\qed 

This gives us another characterization
of harmonic potentials $f$ on $\U\subset\X$.
These are the potentials for which
the associated current is a null divergence
flow on $\U$ (satisfies the first Kirchoff's law)
or, by Stokes Lemma, has zero flux
through any finite cut-set $\partial K$
for which $K\subset \U$.

\medskip\par
We close this section
with a few definitions.
For $A$, $B$ disjoint subsets of $\X$
we say that $\phi$ is a {\em flow from $A$ to $B$}
when
\begin{eqnarray}
\forall a\in A, && \div_a \phi \geq 0\\  
\forall b\in B, && \div_b \phi \leq 0\\  
\forall x\not\in A\cup B, && \div_x \phi = 0
\end{eqnarray}
Any flow $\phi$ is a flow from some $A$ to some $B$.
If $A$ and $B$ are minimal for this property,
we call {\em sources} the elements of $A$
and {\em sinks} those of $B$.
The {\em strength} of a flow $\phi$ with sources in $A$
and sinks in $B$
is
\begin{equation}
|\phi|:= \max\left\{
  \sum_{a\in A} \div_a \phi\: ;\;
  -\sum_{b\in B} \div_b \phi
\right\} 
\end{equation}
A {\em unitary} flow is a flow of strength 1.
If $\phi$ is a unitary flow from $A$ to $B$
and $\X$ is finite then,
by Stokes' lemma with $K=\X$, 
\begin{equation}
\sum_{a\in A}\div_a\phi=
-\sum_{b\in B}\div_b\phi=1
\end{equation}
If $\phi$ is a unitary flow from $A$ to $B$
and $B$ is empty, then we say that $\phi$
is a unitary flow from $A$ to infinity.

\subsection{Equilibrium potential between disjoint subsets}

Consider $A$ and $B$ subsets of $\X$
that satisfy  
\begin{equation}
A\cap B = \emptyset \mbox{ and }\forall x\in\X,\; P_x(\tau_{A\cup B}<+\infty) =1
\label{rc}
\end{equation}
with $P$ the law of the random walk $\xi$
associated with the network.
Assuming that $A$ and $B$
are disjoint subsets of $\X$,
condition~(\ref{rc}) certainly holds when $\xi$ is recurrent,
or
\begin{equation}
\U:=\X\setminus(A\cup B)
\end{equation}
is finite.

Fix the potential at $V_A$ on $A$
and $V_B$ on $B$.
Condition~(\ref{rc}) ensures that Kakutani's
solution of the Dirichlet problem on~$\U$
with such boundary conditions
is well defined. It turns to be
\begin{equation}
V:x\in\X\mapsto V_A P_x(\tau_A<\tau_B) + V_B P_x(\tau_B<\tau_A)
\end{equation}
This is the only one {\em bounded} solution of the Dirichlet
problem.
Indeed, writing~$\X$ as the union
of an increasing sequence of finite sets $K_n$
we can define, for any bounded solution $f$
and all $x\in\X$
\begin{equation}
f_n(x):=E_x\left[f\Big(\xi(\tau_{A\cup B\cup K_n^c})\Big)\right]
\end{equation}
The function $f$ and $f_n$ coincide on $K_n^c$, $A$ and $B$.
Since  both are solutions of a same Dirichlet problem on the finite set $\U\cap K_n$,
they coincide on the whole $\X$.
Now, by dominated convergence, we have
\begin{equation}
f=\lim_{n\rightarrow +\infty} f_n = V
\end{equation}
As a consequence we will refer to $V$ as {\em the} equilibrium potential
conditioned to~$V_A$ on $A$ and $V_B$ on $B$.
In the special case $V_A=1$ and $V_B=0$
we will denote it by $V_{A,B}$:
\begin{equation}
V_{A,B}:x\in\X\mapsto P_x(\tau_A<\tau_B)
\label{dfnvab}
\end{equation}

The current associated with $V$ is (Ohm's law)
\begin{equation}
i = -c\nabla V = -(V_A-V_B)c\nabla V_{A,B}
\end{equation}
its divergence
is zero outside $A\cup B$ (Kirchoff's law),
while for $a$ in $A$ we have
\begin{eqnarray}
\div_a i
&=& -\mu(a)\L_aV\\
&=&(V_A-V_B)\mu(a)(-\L_aV_{A,B})\\
&=&(V_A-V_B)\mu(a)\sum_{y\in \X} p(a,y)[P_a(\tau_A<\tau_B) - P_y(\tau_A<\tau_B)]\\
&=&(V_A-V_B)\mu(a)\sum_{y\in \X} p(a,y)[1- P_y(\tau_A<\tau_B)]\\
&=&(V_A-V_B)\mu(a)\sum_{y\in \X} p(a,y)P_y(\tau_A>\tau_B)\\
&=&(V_A-V_B)\mu(a)P_a(\tau^+_A>\tau^+_B)
\end{eqnarray}
with, for any $S\subset\X$,
\begin{equation}
\tau_S^+:= \min\left\{n>0 :\: \xi(n) \in S\right\}
\end{equation}
The same computation gives for any $b$ in $B$
\begin{eqnarray}
\div_b i
&=&(V_B-V_A)\mu(b)P_b(\tau^+_B>\tau^+_A)
\end{eqnarray}  

By reversibility we have
\begin{eqnarray}
&&\sum_{a\in A} \mu(a)P_a(\tau_A^+ > \tau_B^+)\\
&&\quad = \quad 
\sum_{a\in A}\sum_{b\in B}\sum_{n>0} \mu(a)P_a(\tau_A^+ > \tau_B^+ = n, \xi(n)= b)\\
&&\quad = \quad 
\sum_{a\in A}\sum_{b\in B}\sum_{n>0} \mu(b)P_b(\tau_B^+ > \tau_A^+ = n, \xi(n)= a)\\
&&\quad = \quad 
\sum_{b\in B}\mu(b)P_b(\tau_B^+ > \tau_A^+)
\end{eqnarray}
As a consequence $i$ is a flow of strength
\begin{equation}
|i|= |V_A - V_B|\sum_{a\in A} \mu(a)P_a(\tau_A^+ > \tau_B^+)
=|V_A - V_B|\sum_{b\in B}\mu(b)P_b(\tau_B^+ > \tau_A^+)
\end{equation}

We call {\em capacity} of the pair $(A,B)$ and denote by $C_{A,B}$
the strength of the current associated with $V_{A,B}$
\begin{equation}
C_{A,B}:=\sum_{a\in A} \mu(a)P_a(\tau_A^+ > \tau_B^+)
=\sum_{b\in B}\mu(b)P_b(\tau_B^+ > \tau_A^+)
\label{dfncp}
\end{equation}

Assuming that $C_{A,B}$ is finite, for example when
$A$ or $B$ are finite,
\begin{equation}
i_{A,B}:= \frac{-c\nabla V_{A,B}}{C_{A,B}}
\end{equation}
is a unitary flow from $A$ to $B$.

Writing $\X$ as the union
of an increasing sequence of finite set $K_n$,
replacing $A$ by
$A_n=A\cap K_n$, $B$ by $B_n=K_n^c$
and sending $n$ to infinity
we get an extension of these notions
that turns to be useful when dealing, for example,
with recurrence and transience problems
(Section~\ref{appl1}).
When $n$ goes to infinity $V_{A_n,B_n}$ increases to the limit 
\begin{equation}
h_A:= P(\tau_A < +\infty)
\end{equation}
$C_{A_n,B_n}$ decreases to a non-negative limit,
called {\em capacity of $A$}
\begin{equation}
C_A:= \sum_{a\in A} \mu(a)P_a(\tau_A^+ = +\infty)
\end{equation}
and, if $C_A\in ]0;+\infty[$, then $i_{A,B}$ converges
to a unitary flow $\phi_A$ from $A$ to infinity. 

\section{Energy dissipated in a finite network}

\subsection{Conductance and potentials}

The energy dissipated per time unit in a finite or infinite electrical
network $(\X,c)$ by a potential $f$, or its associated current $i$,
is
\begin{eqnarray}
\D(f) &:=& \frac{1}{2}\sum_{e\in\E} r(e) i^2(e)\\
&=& \frac{1}{2}\sum_{x,y\in\X} c(x,y) [f(x)-f(y)]^2
\end{eqnarray}
The factor $1/2$ is here to ensure
that each pair of connected distinct nodes is counted just once.
$\D(.)$ is the quadratic form associated with the bilinear
{\em Dirichlet form} $\D(.,\!.)$.
As sum of non-negative numbers,
$D(f)$ is always well defined,
even though not always finite.
But the same will not be true
for some of the sums we will write.
To ensure the validity of our next 
calculations we will assume
in this section and the next one that $\X$
is finite.

If $a$ an $b$ under potential $1$ and $0$
are two single
points of an electrical network
made of these two points only,
the energy dissipated in the network
under this potential would be
\begin{equation}
c(a,b)(1-0)^2 = c(a,b)
\end{equation}
This suggests:
\begin{dfn}[Effective conductance]
\label{dfnc}
  If $A$ and $B$ are two disjoint subsets of 
  a finite network $\X$,
  the {\em effective conductance between $A$
  and $B$} is 
  \begin{equation}
  C(A,B) := \D(V_{A,B})
  \end{equation}
($V_{A,B}$ defined in ~(\ref{dfnvab})).
\end{dfn}

If $\X$ were restricted 
to the simple disjoint union
$A\cup B$ each edge of the cutset
$\partial A = - \partial B$
would feel a difference of potential
equal to $1$ and together they would carry
a flow of strength $C_{A,B}$.
As a consequence we would have
\begin{equation}
C(A,B) = \D(V_{A,B}) = C_{A,B}
\end{equation}
This is a general fact:
\begin{prp}
  Capacity and effective conductance coincide.
\end{prp}

\noindent
{\bf Proof:}
Recalling that the current $i$ associated with $V_{A,B}$
is $C_{A,B}.i_{A,B}$ and that $i_{A,B}$ is a unitary
flow from $A$ to $B$, we have:
\begin{eqnarray}
C(A,B) 
& = & \frac{1}{2}\sum_{x,y} c(x,y) [V_{A,B}(x) -V_{A,B}(y)]^2\\
& = & \frac{1}{2}\sum_{x,y} i(x,y)[V_{A,B}(x) -V_{A,B}(y)]\\
& = & C_{A,B}
\sum_{x,y} 
  i_{A,B}(x,y)
  V_{A,B}(x)
\\
& = & C_{A,B}
\sum_{x} 
\underbrace{
  V_{A,B}(x)
}_{
  \stackrel{\scriptstyle 1 \mbox{\small\ on } A,}
  {\scriptstyle 0 \mbox{\small\ on } B}
}
\;
\underbrace{
  \div_x i_{A,B}
}_{
  {\scriptstyle 0 \mbox{\small\ on } A\cup B}
}\\ 
&=& C_{A,B}\sum_{x\in A} \div_x(i_{A,B})\\
&=& C_{A,B}
\end{eqnarray}
\qed

\medskip\par
Effective conductance satisfies
a variational principle:

\begin{prp}[Dirichlet's principle]
\label{varpc}
  \begin{equation}
  C(A,B) = \min\left\{ 
    \D(f) :\:
    f|_A\equiv 1, f|_B\equiv 0
  \right\}
  \end{equation}
  and this minimum is reached in 
  $V_{A,B}$ only.
\end{prp}

\noindent
{\bf Proof:}
Any potential $f$ that is equal to 1
on $A$ and 0 on $B$ can be written
in the form
\begin{equation}
f=V+h
\end{equation}
with
\begin{equation}
V=V_{A,B}, \quad  h|_A\equiv 0, \quad h|_B\equiv 0
\end{equation}
Now 
\begin{equation}
\D(f)=\D(V+h)=\D(V)+\D(h)+2\D(V,h)
\end{equation}
and, denoting by $i$ the current associated with $V$,
\begin{eqnarray}
\D(V,h) &=& \frac{1}{2}\sum_{x,y} c(x,y)[V(x)-V(y)][h(x)-h(y)]\\
&=& \frac{1}{2}\sum_{x,y}i(x,y)[h(x)-h(y)]\\
&=& \sum_{x,y}i(x,y)h(x)\\
&=& \sum_{x}h(x)\div_x i
\end{eqnarray}
Since $h$ equals 0 on $A\cup B$ and $i$ has a null divergence
outside $A\cup B$ we get
\begin{equation}
\D(f)=\D(V)+\D(h)>D(V)
\end{equation} 
as soon as $h\not\equiv 0$.
\qed

\medskip\par
From Dirichlet's principle one gets immediately:
\begin{prp}[Rayleigh's monotonicity law]
  If $c_1\leq c_2$ are such that $(\X,c_1)$ and $(\X,c_2)$ are two finite electrical networks,
  then, for any $A$ and $B$ disjoint subsets of $\X$, 
  $C_1(A,B)\leq C_2(A,B)$, with obvious notation.
\end{prp}
   
We postpone to sections~\ref{appl1}, \ref{appl2}, \ref{appl3}
examples and applications.

\subsection{Resistance and flows}

The energy dissipated per time unit in a finite or infinite
electrical network $(\X,c)$ by a flow~$\phi$
is
\begin{equation}
\D(\phi) := \frac{1}{2}\sum_{e\in\E} r(e) \phi^2(e)
\label{dfndphi}
\end{equation}
If $\phi$ is the current associated
with some potential $f$,
we have, of course,
\begin{equation}
\D(\phi) = \D(f)
\label{co}
\end{equation}
Not all the flows can be derived
from a potential and~(\ref{dfndphi})
generalizes~(\ref{co}).

Consider now $A$ and $B$ two disjoint subsets
of a finite network $\X$.
By the previous variational principle,
any potential that is equal to 1 on $A$ and 0 on $B$
gives an upper bound on $C(A,B)$. 
We derive now a second variational principle
for which any unitary flow
from $A$ to $B$
will give a {\em lower} bound on $C(A,B)$.

\begin{dfn}[Effective resistance]
\label{dfnr}
  If $A$ and $B$ are disjoint subsets of a finite network $\X$, 
  the {\em effective resistance between $A$ and $B$} is 
  \begin{equation}
  R(A,B) := \frac{1}{C(A,B)}
  \end{equation}
\end{dfn}

Effective resistance satisfies
the following variational principle,
cited from~\cite{TT}
by Doyle and Snell \cite{DS}:

\begin{prp}[Thomson's principle]
\label{varpr}
  \begin{equation}
  R(A,B) = \min\left\{
    \D(\phi) :\:
    \mbox{$\phi$ unitary flow from $A$ to $B$}
  \right\} 
  \end{equation}
  and this minimum is reached in $i_{A,B}$
  only.
\end{prp}

\noindent
{\bf Proof:}
The unitary flow $i=i_{A,B}$
is the current associated 
with the potential
\begin{equation}
V=\frac{V_{A,B}}{C(A,B)}
\end{equation}
By bilinearity of the Dirichlet form,
\begin{equation}
\D(i)= \D(V) = \frac{C(A,B)}{C(A,B)^2} = R(A,B)
\end{equation}

Now, any unitary flow from $A$ to $B$, $\phi$,
can be written
\begin{equation}
\phi = i + \delta
\end{equation}
with $\delta$ a flow that satisfies
\begin{eqnarray}
\sum_{a\in A} \div_a \delta &=& 0\\
\sum_{b\in B} \div_b \delta &=& 0\\
\forall x\not\in A\cup B,\;\div_x \delta &=& 0
\end{eqnarray}
so that
\begin{eqnarray}
\D(\phi)
&=& \D(i)+\D(\delta)+\frac{1}{2}\sum_{e\in\E}2r(e)i(e)\delta(e)\\
&=& \D(i)+\D(\delta)+\sum_{(x,y)\in\E}[V(x)-V(y)]\delta(x,y)\\
&=& \D(i)+\D(\delta)+2\sum_{(x,y)\in\E}V(x)\delta(x,y)\\
&=& \D(i)+\D(\delta)+2\sum_{x}V(x)\div_x\delta\\
&=& \D(i)+\D(\delta)+2\sum_{a\in A}C(A,B)^{-1}\div_a\delta\\
&=& \D(i)+\D(\delta)\\
&>& \D(i)
\end{eqnarray}
as soon as $\delta\not\equiv 0$.
\qed

As a consequence, any unitary flow from $A$ to $B$
will give an upper bound on the resistance,
that is a {\em lower} bound on the conductance.
See sections~\ref{appl1}, \ref{appl2}, \ref{appl3} for applications.

\section{Condensers}

\subsection{Capacity and charge} \label{cnc}

Let us go back for a while
to the continuum.
A condenser can be modelized
as a bounded connected open domain $\U$ in $\R^3$
(the domain of the dielectric)
that separates, and is bordered by,
two conductors $A$ and $B$,
at potential $V_A$ and $V_B$.
There cannot be any charge
outside the conductors
and we have
\begin{eqnarray}
\vec{E}&=&-\nabla V\\
\div\vec{E} &=& \frac{\rho}{\epsilon}
\label{maxx}
\end{eqnarray}
Since $V$ is constant on $A$ and $B$
the equations imply that there
cannot be any {\em volumic} charge density.
Physicists say that there can
only be a {\em superficial density}
of charge (on $\partial A$ and $\partial B$)
and using Gauss theorem on an infinitesimal volume
around $a$ in $\partial A$ they conclude
that the superficial density of charge in $a$ is given by
\begin{equation}
q_a = \epsilon(\vec{E}.\vec{n})(a)
\end{equation}
where $\vec{n}$ is the unitary vector
orthogonal to $\partial A$
and directed towards ${\cal U}$.
The total charge on $A$
is given by 
\begin{equation}
Q_A=\int_{\partial A}q_a\;d\sigma(a)
\end{equation}
and the same computation can be reproduced for $B$.
Since potential is defined up to an additive constant
we can replace $V_B$ by 0 and $V_A$ by $V_B-V_A$,
then by linearity of the Dirichlet problem
we get that $Q_A$ depends linearly on $V_A - V_B$, i.e.,
there is a constant $C$, that depends on $A$ and $B$ 
such that
\begin{equation}
Q_A=C (V_A-V_B)
\end{equation}
This constant is called {\em capacity} of the condenser.
In addition the energy contained
in the condenser is given by
\begin{equation}
\int_{\partial A}(V_A-V_B)q_a\, d\sigma(a)
+ \int_{\partial B}0.q_b\, d\sigma(b)
= C (V_A-V_B)^2
\end{equation}

In the context of our electrical network
with $A$ and $B$ that satisfy~(\ref{rc})
under potential $V_A$ and $V_B$ respectively,
for which we know
the equilibrium potential $V$
\begin{equation}
V=V_B  + (V_A-V_B)P_.(\tau_A<\tau_B)
\label{dfnv}
\end{equation}
and the energy dissipated in the network
per time unit 
\begin{equation}
\D(V)= C(A,B)(V_A-V_B)^2
\end{equation} 
the previous considerations lead us:
\begin{itemize}
\item[i)] to define the {\em charge} in any $x\in\X$ 
by analogy with (\ref{maxx}):
\begin{equation}
q_x:=\div_x i
\label{dfnq}
\end{equation}
with $i$ the current associated with $V$.
This is equal to 0 outside $A\cup B$, and for 
$a\in A$, $b\in B$ we get
\begin{eqnarray}
q_a 
&=&(V_A-V_B)\mu(a)P_a(\tau^+_A>\tau^+_B)\\
q_b 
&=&(V_B-V_A)\mu(b)P_b(\tau^+_B>\tau^+_A)
\end{eqnarray}
We recover the ``point-effect'': the higher the escape probability,
the higher the charge.
\item[ii)]
to identify, assuming that $V_A\geq V_B$, the strength
of the current $i$ with the {total charge in $A$}
\begin{equation}
Q_A:=\sum_{a\in A}(V_A-V_B)\mu(a)P_a(\tau^+_A>\tau^+_B)
\end{equation}
and to observe that the two notions of capacity, like
the two notions of energy (contained in the condenser
and dissipated per unit time in the network), coincide in their
probabilistic interpretation, when  any dimensional consideration
disappears. 
\end{itemize}

\medskip\par
We close this section with the
\begin{dfn}[Harmonic measure]
  Given $A$ and $B$ subset of $\X$
  that satisfy~(\ref{rc}) and such that
  \begin{equation}
  C_{A,B}<+\infty
  \end{equation}
  the {\em harmonic measure} on $A$
  is the normalized charge density on~$A$
  under the equilibrium potential $V_{A,B}$
  (or $V$ defined in~(\ref{dfnv})).
  This is the probability measure $\nu_A$ on $A$
  such that, for all $x$ in $A$,
  \begin{eqnarray}
  \nu_A(x) 
  &=& \frac{\mu(x)P_x(\tau_A^+>\tau_B^+)}
  {\sum_{a\in A}\mu(a)P_a(\tau_A^+>\tau_B^+)}\\
  &=& \frac{\mu(x)P_x(\tau_A^+>\tau_B^+)}{C_{A,B}}
  \end{eqnarray}
\end{dfn}

The harmonic measure can be obtained by conditioning the stationary
measure by~$A$ and the
event ``the process $\xi$ 
that start from the sampled point $x$ in $A$
stays outside $A$ at all
positive times before $\tau_B$''.
This does not mean that the random
walk that starts under $\nu_A$ cannot
visit $A$ many times before reaching $B$.
Indeed the conditioning is on the starting
point only: it just means that $\nu_A$
selects the points with the higher
escape probability, not that once chosen
the starting point the escape will occur.
However it is important to note
that $\nu_A$ is concentrated on the internal
border of $A$. 

\subsection{The Green function}

$(\X,c)$ is an electrical network
associated with the Markov chain $\xi$.

\begin{dfn}[Green function]
For any $B\subset\X$ we define the {\em Green function}
\begin{eqnarray}
G_B: (x,y)\in\X^2&\mapsto& E_x\left[\sum_{n=0}^{\tau_B-1} \one_{\{y\}}(\xi(n))\right]\\
&&= \sum_{n\geq 0}P_x\left(\xi(n)=x \mbox{ and } n<\tau_B\right) \label{gsum}
\end{eqnarray}
\end{dfn}

$G_B(x,y)$ is the expected number of visits in $y$
starting from $x$ and before hitting $B$.
Using~(\ref{gsum}) and the reversibility
of $\xi$ we have, for all $x$ and $y$ in $\X$
\begin{equation}
\mu(x)G_B(x,y) = G_B(y,x)\mu(y)
\label{grev}
\end{equation}

If $A$ and $B$ subsets of $\X$ satisfy condition~(\ref{rc})
then the Green function $G_B$
is intimately linked to the potential
$P.(\tau_A<\tau_B)$. To see that we use 
the so-called {\em last exit decomposition}.
We define
\begin{equation}
L_{A,B}:=\sup\left\{
  n \geq 0 :\: \xi(n)\in A \mbox{ and } n<\tau_B
\right\}
\end{equation}
with the usual convention
\begin{equation}
\sup \emptyset = -\infty
\end{equation}
and we have, for all $x$ in $\X$,
using the Markov property and~(\ref{gsum}):
\begin{eqnarray}
P_x(\tau_A<\tau_B)
&=& P_x\left( L_{A,B} \geq 0 \right)\\
&=& \sum_{n\geq 0} P_x\left( L_{A,B}=n \right)\\
&=& \sum_{n\geq 0} \sum_{a\in A} 
P_x\left( \xi(n)=a, n<\tau_B \right)
P_a\left( \tau_A^+>\tau_B^+ \right)\\
&=& \sum_{a\in A} 
G_B(x,a)
P_a\left( \tau_A^+>\tau_B^+ \right)
\end{eqnarray}
that is, by~(\ref{grev}),
\begin{eqnarray}
P_x(\tau_A<\tau_B)
&=& \sum_{a\in A} 
\frac{G_B(a,x)}{\mu(x)}
\underbrace{
  \mu(a)P_a\left( \tau_A^+>\tau_B^+ \right)
}_{
  \mbox{\footnotesize charge in $a$ under $V_{A,B}$}
}\label{vsum}
\end{eqnarray}
In the electrostatic language
we would have say that each charge $q_a$
creates the potential 
\begin{equation}
V^a=\frac{G_B(a,.)}{\mu(.)}q_a
\end{equation}
Indeed, the previous calculation
made in the special case $A=\{a\}$
gives
\begin{equation} 
\frac{G_B(a,.)}{\mu(.)}=\frac{P.\left(\tau_{a}<\tau_B\right)}
{\mu(a) P_a\left(\tau_{a}^+>\tau_B^+\right)}
\end{equation}
so that $V^a$ is harmonic
on $B^c\setminus\{a\}$ (satisfies the local m.v.p.).

Assuming that $C_{A,B}$ is finite, formula~(\ref{vsum}) also gives much
information on the random walk that starts
under the harmonic measure $\nu_A$ and stops
in $B$.
First, it links potential, capacity and stationary measure
with the expected number of visits
to any point $x$ before $\tau_B$.
Multiplying by $\mu(x)$ and dividing
by $C(A,B)$ we get
\begin{equation}
E_{\nu_A}\left[\sum_{n<\tau_B}\one_{\{x\}}(\xi(n))\right]
= \frac{\mu(x)P_x(\tau_A<\tau_B)}{C_{A,B}}
\label{nbvis}
\end{equation}

Second, summing over all $x$ outside $B$,
we get the expected hitting time of $B$.
This is the main formula
that was introduced in \cite{BEGK1}
for the study of metastability:
\begin{equation}
  E_{\nu_A}\left[\tau_B\right] =
  \frac{1}{C_{A,B}}
  \sum_{x\not\in B}\mu(x)P_x(\tau_A<\tau_B)
  = \frac{\mu(V_{A,B})}{C_{A,B}}
\label{xmeta}
\end{equation}

Last, it makes possible
to give the probabilistic interpretation
of the unitary flow $i_{A,B}$.
For $e=(x,y)\in\E$ we have
\begin{eqnarray}
i_{A,B}(e)
&=& c(x,y)\left[\frac{P_x(\tau_A<\tau_B)}{C_{A,B}}-\frac{P_y(\tau_A<\tau_B)}{C_{A,B}}\right]\\
&=& \sum_{a\in A}\nu_A(a)\left(
  \frac{G(a,x)}{\mu(x)}
  - \frac{G(a,y)}{\mu(y)}
\right)c(x,y)\\
&=& \sum_{a\in A}\nu_A(a)\Big(
  {G(a,x)}{p(x,y)}
  - {G(a,y)}{p(y,x)}
\Big)\\
&=&
E_{\nu_A}\left[\sum_{n<\tau_B}\left(\one_{\{e\}}-\one_{\{-e\}}\right)(\xi(n),\xi(n+1))\right]
\end{eqnarray}
This is the expected net flux 
of the walk through $e$. 

\section{Application to transience and recurrence} \label{appl1}

\subsection{Recurrence and conductance} 

Let $\xi$ be a reversible ergodic Markov chain on $\X$,
and $(\X,c)$ an associated electrical network.
The random walk is {\em recurrent} if
\begin{equation}
\forall x,y \in \X,\; P_x(\tau_y<+\infty)=1
\end{equation}
otherwise it is {\em transient}.
If $\X$ is finite $\xi$ (that is assumed
to be ergodic) is necessarily recurrent.
In general we can write $\X$ as union
of an increasing sequence of finite connected subsets $K_n$,
and we have

\begin{prp}\label{crit}
  The following assertions are equivalent:
  \begin{eqnarray}
  i) &&\mbox{$\xi$ is recurrent}\\ 
  ii) && \exists a \in \X,\; P_a\left(\tau_a^+ < +\infty\right)=1\\
  iii) && \exists a \in \X,\; 
  E_a\left[{\textstyle \sum_{n\geq 0}} \one_{\{a\}}\Big(\xi(n)\Big)\right] = +\infty\\
  iv) && \exists a\in \X, \; \lim_{n\rightarrow +\infty} G_{K_n^c}(a,a) = +\infty\\
  v) && \exists a\in \X, \; \lim_{n\rightarrow +\infty} C_{a,K_n^c} = 0\\
  vi) && \exists a\in \X, \; \exists n_0\geq 0,\; 
         \exists (f_n)_{n \geq n_0}: \X \rightarrow [0;1],\;\nonumber\\
     &&\qquad
     f_n(a) = 1,\; f_n|_{K_n^c}\equiv 0,\; \lim_{n\rightarrow +\infty}\D(f_n) = 0 
  \end{eqnarray}
\end{prp}
 
\noindent
{\bf Proof:} i) $\Rightarrow$ ii) is clear and ii) $\Rightarrow$ i), since,
for all $x\in \X$,
\begin{equation}
P_a(\tau_x<+\infty) >0
\end{equation}
and a random walk that almost surely visits $a$ infinitely many times,
will almost surely visit $x$.

ii) $\Rightarrow$ iii) is clear and the number of visits in $a$
for the random walk that starts in $a$ is distributed like
a geometric variable of parameter
\begin{equation}
p=P_a(\tau_a^+=+\infty)
\end{equation}
If $p\neq 0$ the expected number of visits in $a$ is finite
and iii) $\Rightarrow$ ii) follows.

We have iii) $\Leftrightarrow$ iv) by Beppo Levi's theorem,
and iv) $\Leftrightarrow$ v) follows from~(\ref{nbvis}) applied
with $A=\{a\}$ and $B = K_n^c$.

v)$\Leftrightarrow$vi)
follows from a variational principle 
for effective conductances.
It was proved (Proposition~\ref{varpc})
for finite networks and we can extend
it to our situation:
for any $n$ we build a finite network $(\X_n,c_n)$
by collapsing in a single point $b$,
all the nodes in $K_n^c$:
we define $\X_n$ as the union of $K_n$
with a singleton $\{b\}$
and we define $c_n$ by 
\begin{equation}
c_n(x,y):=\left\{
  \begin{array}{ll}
    0 & \mbox{if $(x,y)=(b,b)$}\\
    c(x,y) & \mbox{if $(x,y)\in K_n\times K_n$}\\
    {\sum_{y\in K_n^c} c(x,y)}
    &\mbox{if $(x,y)\in K_n\times \{b\} \cup \{b\}\times K_n$}
  \end{array}
\right.
\end{equation}
On the one hand the law of the random walks 
$\xi_n$, associated
with $(\X_n,c_n)$,
and $\xi$ that start in $a$
are the same up to $\tau_{K_n^c}$.
The total weight in $a$
is the same in the two networks,
hence the capacity
$C_{a,K_n^c}$ associated with $(\X,c)$,
coincides with $C_{a,b}$.
On the other hand $C(a,b)$ satisfies
the variational principle of Proposition~\ref{varpc},
and the Dirichlet form of the corresponding test functions
coincide with the Dirichlet form of the test functions
of what would be the analogous variational principle on $\X$.
This proves the validity of this variational
principle and concludes the proof.
\qed

\medskip\par\noindent
{\bf Example:}
For the simple random walk on $\Z^2$
the conductance of each edge is $1/4$. 
We set, for all $n\geq 1$,
\begin{equation}
K_n := [-(n-1);n-1]^2
\end{equation}
and we consider the potentials
\begin{equation}
\begin{array}{rcrcl}
  f_n  & : & \Z^2 & \longrightarrow & [0;1]\\
  &&x&\longmapsto&\left\{
    \begin{array}{ll}
      1 - \frac{\ln(1+\|x\|_\infty)}{\ln(1+n)} &\mbox{if $x \in K_n$}\\
      0 &\mbox{if $x \in K_n^c$}
    \end{array}
  \right.
\end{array}
\end{equation}
We have
\begin{eqnarray}
\D(f_n)
&=& \frac{1}{\ln^2(1+n)}\sum_{k = 1}^{n}
\frac{(8k-8)\vee 1}{4}
\left[
   \ln(k+1) -\ln k
\right]^2\\
&\leq& \frac{1}{\ln^2(1+n)}\sum_{k = 1}^{n}
2k\frac{1}{k^2}\\
&\leq& 2\frac{1+\ln(n+1)}{\ln^2(1+n)}
\end{eqnarray}
and we conclude that the random walk is recurrent.

This may not be the simplest proof 
of the recurrence,
but it is the most resistant I know.
For example if we remove any set of edges
from the initial graph,
then by Rayleigh's monotonicity law
the random walk obtained by refusing
the jump each time it tries to move
along a removed edge is recurrent
on each connected component of
the obtained graph.

\subsection{Lyons' criterion}

We can add to our list of Proposition~\ref{crit}
a last criterion, due to T.~Lyons (see~\cite{Ly}),
for deciding whether a given random walk $\xi$
is recurrent or not.
\begin{prp}
  A reversible Markov chain $\xi$ associated with
  an electrical network  $(\X,c)$
  is transient
  if and only if there is a unitary flow
  from some $a$ in $\X$ to infinity
  that dissipates a finite energy in the network.
\end{prp}

\noindent
{\bf Proof:}
If $\xi$ is transient then, for any $a$
in $\X$, 
\begin{equation}
\lim_{n\rightarrow +\infty} C_{a,K_n^c} = C_a > 0
\end{equation}
In this case $(i_{a,K_n^c})_{n\geq 0}$ converges
to a unitary flow $\phi_a$ from $a$ to infinity
and we have 
\begin{equation}
\D(\phi_a)=\lim_{n\rightarrow +\infty} \D(i_{a,K_n^c})
=\lim_{n\rightarrow +\infty} \D\left(\frac{V_{a,K_n^c}}{C_{a,K_n^c}}\right)
=\lim_{n\rightarrow +\infty} \frac{C_{a,K_n^c}}{C^2_{a,K_n^c}}
= \frac{1}{C_a} <+\infty
\end{equation}  

If there is a unitary flow $\phi$ 
from $a\in\X$ to infinity with
\begin{equation}
\D(\phi) <+\infty
\end{equation}
then, for $n$ large enough, $\phi$ is also a unitary flow 
from $a$ to any $K_n^c$ and,
denoting by $\D_n$ the Dirichlet form
on the network $(\X_n, c_n)$ built
by collapsing $K_n^c$ in a single point $b$,
by $r_n$ the associating resistances,
and defining a unitary flow from $a$ to $b$ by
\begin{equation}
\phi_n(e):=\left\{
  \begin{array}{ll}
    \phi(e) & \mbox{if $e\in\E\cap K_n\times K_n$}\\
    \sum_{\stackrel{\scriptstyle y\in K_n^c}{\scriptstyle (x,y)\in\E}} \phi(x,y)
    &\mbox{if $e=(x,b)$ with $x\in K_n$}
  \end{array}
\right.
\end{equation}
we have (using Jensen's inequality to get ~(\ref{jensen})):
\begin{eqnarray}
\D(\phi)
&=& \frac{1}{2}\sum_{e\in\E} r(e) \phi^2(e)\\
&\geq& \frac{1}{2}\sum_{e\in\E\cap K_n\times K_n} r(e) \phi^2(e)
       + \sum_{x\in K_n} \sum_{\stackrel{\scriptstyle y\in K_n^c}{\scriptstyle (x,y)\in\E}} r(x,y) \phi^2(x,y)\\
&=& \frac{1}{2}\sum_{e\in\E\cap K_n\times K_n} r(e) \phi^2(e)
    + \sum_{x\in K_n} \sum_{\stackrel{\scriptstyle y\in K_n^c}{\scriptstyle (x,y)\in\E}}
      c(x,y) \left(
        \frac{\phi(x,y)}{c(x,y)}
      \right)^2\\
&\geq& \frac{1}{2}\sum_{e\in\E\cap K_n\times K_n} r(e) \phi^2(e)
       + \sum_{x\in K_n} c_n(x,b)\left(
           \frac{
             \phi_n(x,b)
           }{
             c_n(x,b)
           }
         \right)^2
       \label{jensen}\\
&=& \frac{1}{2}\sum_{e\in\E\cap K_n\times K_n} r(e) \phi^2(e)
    + \sum_{x\in K_n} r_n(x,b)\phi^2(x,b)\\
&=& D_n(\phi_n)\\
&\geq& R(a,b)\\
&=&C^{-1}_{a,K_n^c}
\end{eqnarray}
and we conclude that $C_{a,K_n^c}$ decreases with $n$ towards a strictly positive value. 
\qed

\medskip\par\noindent
{\bf Example:}
Consider the simple random walk on $\Z^d$
with $d\geq 3$. We can build
a unitary flow $\phi$ from 0 to infinity in the following way.
First we associate with each $\theta$ in 
\begin{equation}
S^{d-1} := \partial B_2(0,1)
\end{equation}  
a path $\gamma$ from 0 to infinity, $e^1$, $e^2$,~\dots 
such that $(\|e^k_+\|_2)_{k\geq 0}$ is increasing
and, for all $k\geq 0$, the distance between
$e_+^k$ and the half line $[0,\theta)$ is less
than 2.
Second we define, for all $e\in\E$,
\begin{equation} 
\phi^\theta(e):=\one_\gamma(e) - \one_\gamma(-e)
\end{equation}
$\phi^\theta$ is a unitary flow from 0 to infinity.
Last we define $\phi(e)$ as the expected
value of $\phi^\theta(e)$ when $\theta$ is chosen
according to the uniform probability measure $\P$ on $S^{d-1}$:
\begin{equation}
\phi(e):=\P(e\in\gamma)-\P(-e\in\gamma)
\end{equation}
$\phi$ is a unitary flow from 0 to infinity,
we have
\begin{eqnarray}
\D(\phi)
&=& \sum_{e\in\E} r(e)\phi^2(e)\\
&\leq& \cst\sum_{r\geq 1} r^{d-1}\left(\frac{1}{r^{d-1}}\right)^2\\
&=& \cst\sum_{r\geq 1} \frac{1}{r^{d-1}}\\
&<&+\infty
\end{eqnarray}
and we get the transience of the random walk.

Lyons' criterion for transience has proven
to be extremely powerful.
It has been used for example in~\cite{GKZ}
to prove the transience of the random walk
on the infinite supercritical percolation cluster
in dimension $d\geq 3$.

\section{Application to metastability} \label{appl2}

\subsection{Restricted ensemble}

Metastability is characterized by (at least)
two different time scales, a short and a long one,
and an apparent equilibrium.
If the equilibrium of the system is described
by a measure $\mu$,
this apparent equilibrium
is described by a {\em restricted ensemble}
$\mu_{\cal R}$, that is the equilibrium measure
conditioned to a subset ${\cal R}$ of the state space $\X$.
With a probability of order 1,
the system initially described
by a metastable equilibrium $\mu_{\cal R}$
will escape from ${\cal R}$
on the long time scale,
then, on the short time scale,
will go far away from ${\cal R}$
(far away in the sense that 
he will come back to ${\cal R}$
on a third and still longer time scale)
towards a more stable equilibrium.

Such a behaviour can be modelized
through that of an ergodic continuous time Markov process
$X$ on a finite state space $\X$
on which are defined an Hamiltonian $H$ and its associated
Gibbs measure $\mu$ at inverse temperature $\beta>0$
\begin{equation}
\mu:=\frac{1}{Z}\exp\{-\beta H\}
\quad
\mbox{with}
\quad
Z:=\sum_{x\in\X} \exp\{-\beta H(x)\}
\end{equation}
and for which
$X$ is reversible with respect to $\mu$
(so that $\mu$ is the unique equilibrium measure).
The previous expressions ``short and long time scales'',
``probability of order 1''
make then sense in some asymptotic regime, for example
when $\beta$, $|\X|$ or some other parameter
of the dynamic goes to infinity.

In what follows we will consider continuous time Markov processes $X$ defined
by a Metropolis algorithm associated
with $H$, i.e., with a generator of the form
\begin{equation}
L_x f := \sum_{y\sim x} \exp\left\{-\beta\left[H(y)-H(x)\right]_+\right\} (f(y)-f(x))
\label{dfnmet}\end{equation}
(note that~(\ref{dfnmet}) guarantees the reversibility with respect to $\mu$)     
and we will consider the (by far easier) regime $\beta\rightarrow +\infty$,
or a joint regime in which $\beta$ {\em and} $|\X|$ go to infinity. 
We will refer to these two kinds of regimes as {\em finite} and {\em large volume dynamics}
respectively.

\subsection{Finite volume dynamics}
\label{fv}

Our two main examples are Glauber and Local Kawasaki dynamics.
Given $\Lambda$ a finite square box in $\Z^d$
with $d\geq 2$ the Glauber dynamics is defined on the state space
\begin{equation}
\X = \{-1;+1\}^\Lambda
\end{equation}
with Ising Hamiltonian with periodic boundary conditions
\begin{equation}
H:\sigma\in\X\mapsto
-\frac{1}{2}\sum_{\stackrel{\scriptstyle \{i,j\}\subset\Lambda}{\scriptstyle d^T_1(i,j) =1}}J\sigma_i\sigma_j 
-\frac{1}{2}\sum_{i\in\Lambda} h\sigma_i
\end{equation}
where $J>0$ is the ferromagnetic interaction constant, $h>0$ the magnetic field,
and $d^T_1(i,j)$ gives the 1-distance on the torus between the projections of $i$
and $j$.  
It is a single spin flip dynamic,
that is $y\sim x$ in~(\ref{dfnmet})
means that $y$ is obtained from $x=\sigma\in\X$
by changing the value of $\sigma$ in one site $i$ of the torus.

The Local Kawasaki dynamics is defined on the state space
\begin{equation}
\X = \{0;1\}^\Lambda
\end{equation}
with Hamiltonian 
\begin{equation}
H:\eta\in\X\mapsto
\sum_{\stackrel{\scriptstyle \{i,j\}\subset\Lambda\setminus\partial_-\Lambda}{\scriptstyle d_1(i,j) =1}}-U\eta_i\eta_j 
+\sum_{i\in\Lambda} \Delta\eta_i
\end{equation}
where $-U<0$ is the binding energy and $\Delta>0$ an activity parameter.
It is a (locally conservative) nearest neighbours exchange dynamic
with creation and annihilation of particles on the internal border of the box,
that is $y\sim x$ in~(\ref{dfnmet})
means that $y$ is obtained from $x=\eta\in\X$
by exchanging the value of $\eta$ between two nearest neighbour sites $i$ and $j$ in $\Lambda$
or by changing the value of $\eta$ in one site $i\in\partial_-\Lambda$.

Whatever the model we consider, 
the individuation of a set ${\cal R}$
with the previously described properties
is part of the problem.
For finite volume Glauber dynamics
it was done by Neves and Schonmann 
in \cite{NS1}, \cite{NS2} and this was generalized
to a host of situation including
that of the beautiful paper of Schonmann and Shlosman
\cite{SS} that consider metastability for Glauber
dynamics in infinite volume at finite temperature
and in the regime $h\rightarrow 0$.
For finite volume Local Kawasaki dynamics
it was done by den Hollander, Olivieri and Scoppola
in \cite{dHOS} for $d=2$, by den Hollander, Nardi, Olivieri and Scoppola
in cite \cite{dHNOS} for $d=3$.
Assuming that
\begin{equation}
\frac{2J}{h}\:, \frac{U}{2U-\Delta} \quad\in\quad ]1;+\infty[\:\cap\:\N^c
\end{equation}
and defining the {\em critical length} $l_c$ by
\begin{equation}
l_c:=\left\{
  \begin{array}{cl}
    \left\lceil\frac{2J}{h}\right\rceil &\mbox{for Glauber dynamics}\\
    \\
    \left\lceil\frac{U}{2U-\Delta}\right\rceil & \mbox{for Kawasaki dynamics}
  \end{array}
\right.
\end{equation} 
one can define a {\em gate} ${\cal G}$,
set of critical configurations at a same energy $H({\cal G})$
that, for Glauber dynamics and $d=2$,
are  the quasi-squares droplets of $+1$ 
in $\Lambda\setminus\partial_-\Lambda$
of dimensions $(l_c-1)\times l_c$
with a protuberance attached on the long side,
while, for Local Kawasaki dynamics, have for prototype
the quasi-squares droplets of $1$  of dimensions $(l_c-1)\times l_c$
with a protuberance and an extra free particle.

\begin{figure}[htbp]
\begin{center}
\includegraphics
[width=8cm]
{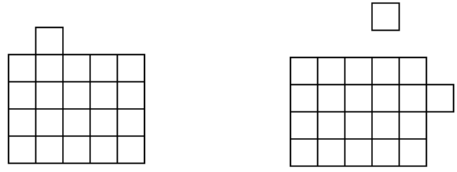}
\caption{Critical configurations for Glauber
and Local Kawasaki dynamics in the case $l_c=5$.}
\end{center}
\end{figure}

\noindent
Then it was shown 
(see in particular \cite{GOS}
for Local Kawasaki)
that, with $a\in\X$ the configuration
made of $-1$ only ($0$ only) and $b\in\X$ the configuration
made of $+1$ only (1 in $\Lambda\setminus\partial_-\Lambda$
and 0 in $\partial_-\Lambda$) for Glauber (Local Kawasaki) dynamics,
for $\Lambda$ large enough:
\begin{description}
\item[{\bf (P1)}] $b$ is the only one fundamental state, that is the global
minimum of the Hamiltonian,
\item[{\bf (P2)}] $a$ is the only one metastable state
in the sense of~\cite{MNOS} that is, with
\begin{equation} 
\Gamma:=\min\left\{
  H(a) \vee \max_{e\in\gamma} H(e_+) :\: 
  \mbox{$\gamma$ is a path from $a$ to $b$}
\right\}-H(a)
\label{minmax}
\end{equation}
there is, for all $x$ in $\X\setminus\{a\}$,
\begin{equation}
\Gamma > \min\left\{
  H(x) \vee \max_{e\in\gamma} H(e_+) :\: 
  \mbox{$\gamma$ is a path from $x$ to $b$}
\right\}-H(x)
\end{equation}
\item[{\bf (P3)}] ${\cal G}$ has the {\em gate property}, i.e.,
any path that realizes the min-max
~(\ref{minmax}) has to cross ${\cal G}$
and reaches its maximum in ${\cal G}$, 
so that, in particular,
\begin{equation}
\Gamma = H({\cal G})-H(a)
\label{gammaa}
\end{equation}
(note that $\Gamma$ depends on $J$ and $h$ or $U$ and $\Delta$ only),
\item[{\bf (P4)}] if $A$ and $B$ are the two cycles 
in the sense of Wentzell and Freidlin~\cite{FW}
that are the connected components of $a$ and $b$
in $H^{-1}(]-\infty,H({\cal G})[)$ then 
\begin{equation}
{\cal G}\subset \partial_+ A
\end{equation}
and, by (P2), 
\begin{equation}
\forall x\not\in B,\;
H(x) > H(a)
\label{cut}
\end{equation}
\end{description}

This is a big amount of information --
(P2) includes a control
of the global energy landscape --
and at this point there are many possible
choices for the set~${\cal R}$.
Natural choices include
\begin{itemize}
\item $A$,
\item  the connected component
of $A$ in $H^{-1}(]-\infty,H({\cal G})])\setminus {\cal G}$,
\item $H^{-1}(]-\infty,H({\cal G})])\setminus {B}$,
\item larger sets (including for example
a small piece of $B$)...
\end{itemize} 
With any of  these choices Wentzell-Freidlin theory
leads, for all $\delta >0$, to 
\begin{equation}
\limsup_{\beta\rightarrow +\infty}
\frac{1}{\beta}\ln P_{\mu_{\cal R}}\left(
  \tau_{{\cal R}^c}, \tau_b 
  \not\in 
  [e^{(\Gamma-\delta)\beta}; e^{(\Gamma+\delta)\beta}]
\right)
<0
\label{metar}
\end{equation}
This does not say much
about the existence of our ``short time scale''
but it is a strong indication
that our ``long time scale''
should be  $e^{\Gamma\beta}$.
It takes, indeed, for the system
initially under $\mu_{\cal R}$, essentially the same (long) time
to reach ${\cal R}^c$ and $b$, and when
the system is in $b$ it is ``far from ${\cal R}$'':
one can see, using reversibility, that, typically,
the system needs at least a time of order
\begin{equation}
e^{(H({\cal G}) -H(b))\beta}
= e^{(\Gamma + H(a) -H(b))\beta}
\gg
e^{\Gamma\beta}
\end{equation}
to go back to ${\cal R}$.

In addition, Wentzell Freidlin theory,
leads also to 
\begin{eqnarray}
\Gamma 
&=& \lim_{\beta\rightarrow +\infty}
    \frac{1}{\beta}\ln E_{\mu_{\cal R}}\left[
      \tau_{{\cal R}^c}
    \right]\label{1mn}\\
&=& \lim_{\beta\rightarrow +\infty}
    \frac{1}{\beta}\ln E_{\mu_{\cal R}}\left[
      \tau_{b}
    \right]\\
&=& \lim_{\beta\rightarrow +\infty}
    \frac{1}{\beta}\ln E_{a}\left[
      \tau_{{\cal R}^c}
    \right]\\
&=& \lim_{\beta\rightarrow +\infty}
    \frac{1}{\beta}\ln E_{a}\left[
      \tau_{{b}}
    \right]\label{4mn}
\end{eqnarray}

We refer to~\cite{MNOS}
for the derivation of all these results
on the basis of (P1)-(P4).

\noindent
{\bf Remarks:}  
{\bf i)} The pathwise approach on which were
based the proofs of (P1)-(P4) gives also
the ``short time scale''.
For example, in the case of the Local Kawasaki dynamics, 
it is of order $e^{(2\Delta-U)\beta}$ \cite{GOS}.
\smallskip\par\noindent
{\bf ii)} Equations~(\ref{1mn})-(\ref{4mn})
are stronger than~(\ref{metar}) in the sense
that the former imply (with Markov inequality)
the upper bound on $\tau_{\cal R}$ and ${\tau_{b}}$
expressed by the latter, while the lower bound
on these times is easy to get using reversibility.
However it is important to note
that~(\ref{1mn})-(\ref{4mn}) imply in general an information 
(of the kind of (P2))
on the global energy landscape. 
If $a$ does not lie on the  bottom
of the deepest well 
(like expressed in (P2))
and can reach, without going in $b$, a well with a depth $\Gamma'$ 
larger than $\Gamma$
with a probability exponentially larger than $e^{-(\Gamma'-\Gamma)\beta}$
then~(\ref{1mn})-(\ref{4mn}) cannot hold.
By contrast, results like~(\ref{metar})
can be derived by a strictly pathwise approach
(see for example \cite{BAC} on Glauber dynamics
in dimension 3)
without such kind of information on the global
energy landscape.

\subsection{Beyond exponential asymptotics}
\label{sharp}

On the basis on (P1)-(P4), potential theory
can improve~(\ref{1mn})-(\ref{4mn})
{\em beyond exponential asymptotics}.
Bovier and Manzo did that in \cite{BM}
and gave the exact asymptotics of 
$E_a[\tau_b]$ for Glauber dynamics.

They did so applying~(\ref{xmeta})
to the sets $\{a\}$ and $\{b\}$,
and, this implied, in particular,
giving some estimates on the capacity.
As a far as the upper bounds (on the capacity)
are concerned,
they estimated $C(a,b)$
with
\begin{equation}
C(a,b)\leq C(A,B)
\end{equation}
where $A$ and $B$ are the two cycles
defined in (P4) (by Dirichlet's principle
the conductance is increasing in its arguments).
To give a lower bound on the capacity
they drop some terms in the Dirichlet form
of the variational principle
(Rayleigh's monotonicity law)
to get a linear network
for which they were able to compute the capacity.
This is equivalent to building
a linear flow and using
Thomson's principle.

We will use a slightly different strategy:
we will apply~(\ref{xmeta})
directly to our cycles $A$ and $B$.
But before doing that we have to pass through
a little
algebra to link the study of our continuous
time Markov process $X$ to that
of the discrete time random walks $\xi$
we dealt with in the previous sections.

Observe that the generator $L$ defined
in~(\ref{dfnmet}) cannot be written
in the form that $\L$ assumed in~(\ref{dfnl}):
given $x\in\X$ the sum on $y\sim x$ of the rates 
\begin{equation}
\lambda(x,y):=\exp\left\{-\beta[H(y)-H(x)]_+\right\}
\end{equation}
is in general larger than one.
But it is certainly smaller than
\begin{equation}
N:=\left\{
\begin{array}{l}
\mbox{number of sites in $\Lambda$ for Glauber}\\
\mbox{number of bonds inside $\Lambda\setminus\partial_-\Lambda$
and sites in $\partial_-\Lambda$ for Kawasaki}
\end{array}
\right.
\end{equation}
We define then the network
$(\X,c)$ with, for all $x$ and $y$ in $\X$,
\begin{equation}
c(x,y):=
\left\{
  \begin{array}{ll}
    0 & \mbox{if $y\neq x$ and $y\not\sim x$}\\
    \frac{\mu(x)\lambda(x,y)}{N} & \mbox{if $y\neq x$ and $y\sim x$}\\   
    \mu(x) - \sum_{\stackrel{\scriptstyle y\sim x}{y\neq x}}
             \frac{\mu(x)\lambda(x,y)}{N} & \mbox{if $y=x$}
  \end{array}
\right.
\label{dfncn}
\end{equation}
Since, for all $x$ in $\X$,
\begin{equation}
  \sum_{y\in\X} c(x,y) = \mu(x)
\end{equation}
the random walk $\xi$ associated with $(\X,c)$
is reversible with respect to $\mu$.
Its generator is defined by 
\begin{equation}   
\L_x f:= \sum_{y\in \X} \frac{c(x,y)}{\mu(x)} (f(y)-f(x))
= \sum_{y\sim x} \frac{\lambda(x,y)}{N} (f(y)-f(x))
= \frac{1}{N}L_x f
\label{rs}
\end{equation}
Recall that we called ``generator of a discrete time Markov chain $\xi$''
that of the continuous time process that updates its position
at each ring of a Poissonian clock of intensity 1
according to the transition probabilities of $\xi$.
Denoting by $\tilde\xi$
this continuous time process~(\ref{rs})
means that $\tilde\xi$ is nothing
but the rescaled process $X$:
$\tilde\xi$ behaves like $X$
except for the fact that it is $N$ times slower.

As a consequence
\begin{equation}
E_{\nu_A}[\tau_B(X)]
= \frac{1}{N}E_{\nu_A}[\tau_B(\tilde\xi)]
= \frac{1}{N}E_{\nu_A}[\tau_B(\xi)]
\end{equation}
and~(\ref{xmeta}) gives
\begin{equation}
E_{\nu_A}[\tau_B(X)]
= \frac{1}{NC(A,B)}\sum_{x\not\in B}\mu(x)P_x(\tau_A<\tau_B)
\label{obj}
\end{equation}
where we have to recall that the conductances
that are involved in the computation of $C(A,B)$
are defined in~(\ref{dfncn}) and depend on $N$ too.

By~(\ref{cut}) 
the last sum in~(\ref{obj}) 
is equivalent to $\mu(a)$
and, in the case of the 2-dimensional
Glauber dynamics, it turns out that
for all $g$ in ${\cal G}$, $A\cup\{g\}\cup B$
is a connected set, while for all distinct
$g$, $g'$ in ${\cal G}$, $g$ and $g'$
are not connected.
These two properties make then
the capacity in~(\ref{obj}) quite
easy to estimate. Indeed, we first note
that for all $x\neq y$  with $x\sim y$,
\begin{equation}
c(x,y)= \frac{\mu(x)\lambda(x,y)}{N}
=\frac{\exp\{-\beta(H(x) \vee H(y))\}}{NZ}
\end{equation}
This implies that,
for any function $f$ on $\X$
that takes its values in $[0,1]$,
is equal to $1$ on $A$ 
and to $0$ on $B$,
all the terms in the Dirichlet from $D(f)$
that involve a node beyond the energy
level $H({\cal G})$
are exponentially smaller than $D(f)$,
and, using the gate property
of ${\cal G}$ and the fact that
for all $g$ in ${\cal G}$, $A\cup\{g\}\cup B$
is a connected set,
we conclude that the $C(A,B)$ 
is equivalent to the capacity of the pair $(A,B)$
in the network $(A \cup {\cal G} \cup B, c)$.
In this network a fraction $2/l_c$
of the nodes in ${\cal G}$ have only one edge towards
$A$ and one edge towards $B$,
while the other nodes have only
one edge towards $A$ and two edges towards $B$.

\begin{figure}[htbp]
\begin{center}
\includegraphics
[width=12cm]
{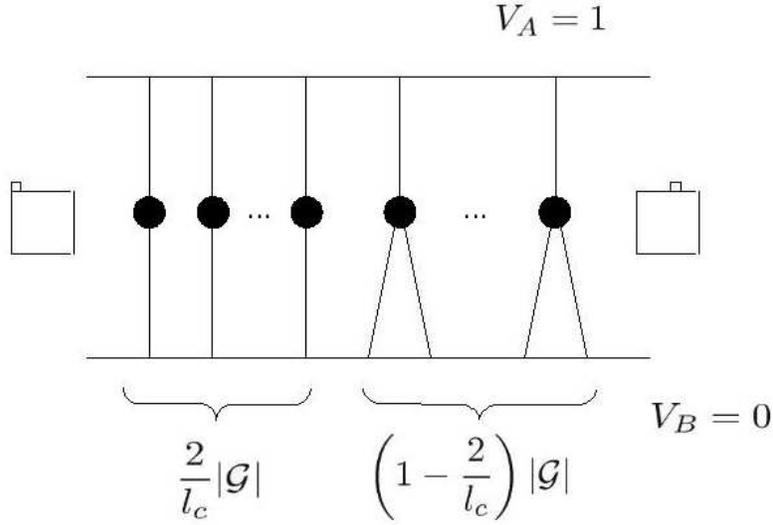}
\end{center}
\caption{The network $(A \cup {\cal G} \cup B, c)$.}
\end{figure}

\noindent
All these edges have the same conductance 
\begin{equation}
\bar c = \frac{\exp\{-\beta H({\cal G})\}}{NZ}
\end{equation}
and we get 
\begin{equation}
C(A,B)
\sim 
\frac{\frac{2}{l_c}|{\cal G}|}{\frac{1}{\bar c} +\frac{1}{\bar c}}
+ \frac{\left(1-\frac{2}{l_c}\right)|{\cal G}|}{\frac{1}{\bar c} +\frac{1}{2\bar c}}
=\frac
{
(2l_c-1)|{\cal G}|
{\exp\{-\beta H({\cal G})\}}
}
{3l_c{NZ}}
\end{equation}
We conclude, using~(\ref{gammaa}),
\begin{equation}
E_{\nu_A}[\tau_B(X)]\sim\frac{{3l_c{Z}}\mu(a)}{{(2l_c-1)|{\cal G}|{\exp\{-\beta H({\cal G})\}}}}
=\frac{{3l_ce^{\Gamma\beta}}}{(2l_c-1)|{\cal G}|}
\end{equation}
From this we get the same estimate on $E_a[\tau_b]$. Here is the logic
of the argument. The probability $\nu_A$ is concentrated on $A$
that is a cycle of depth $\Gamma$,
as a consequence the system will typically reach $a$
in a time exponentially smaller than $e^{\Gamma\beta}$
before going to $B$. But $B$ is a cycle with internal resistance
smaller than $\Gamma$ ($a$ is the only one metastable state),
hence, after reaching $B$ the system will typically go to $b$
in a time exponentially smaller than $e^{\Gamma\beta}$.
This leads, for any small enough $\delta>0$, to 
\begin{equation}
E_a[\tau_b] + o(e^{(\Gamma-\delta)\beta}) = 
E_{\nu_A}[\tau_B(X)]\sim
\frac{{3l_ce^{\Gamma\beta}}}{(2l_c-1)|{\cal G}|}
\end{equation}
and 
\begin{equation}
E_a[\tau_b] \sim
\frac{{3l_ce^{\Gamma\beta}}}{(2l_c-1)|{\cal G}|}
\label{cl}
\end{equation}
To put the argument properly
you have to quantify the probability of ``atypical behaviours''
to control the expectations,
and this, knowing (P1)-(P4), is elementary classical Wentzell-Freidlin theory.

With~(\ref{cl}) everything boils down to the computation
of the number of critical configurations.
For the 2-dimensional finite volume Glauber dynamics
we find
\begin{equation}
|{\cal G}|= 4 l_c |\Lambda|
\end{equation}
(there are $|\Lambda|$ choices
for the south-west corner of the quasi-square,
2 choices for its orientation and
$2l_c$ choices for the position of the protuberance).
We refer to~\cite{BM} for the study of the dynamics
in higher dimension.

For Local Kawasaki the situation is more complex:
first ${\cal G}$ (that is not uniquely defined)
is not so simple, second the electrical network
that connects $A$ and $B$ is ``stretched''
and much more intricate.
But the same method can be applied,
$C(A,B)$ can be estimated
via our two variational principles,
and, once again, everything is reduced
to some computation of $|{\cal G}|$.
This is the difficult point
of~\cite{BdHN} that gives sharp asymptotics
of $E_a[\tau_b]$ for this model
in dimensions~2 and 3.

\subsection{Large volume dynamics}
\label{metal}

Glauber and Kawasaki dynamics
in large volume are defined
as continuous time Markov chains $X$
on the space $\X$ of the configurations
made on -1 and +1 (Glauber) or 0 and 1
(Kawasaki) on the d-dimensional
discrete torus $\Lambda_\beta$
of volume
\begin{equation}
|\Lambda_\beta| = e^{\Theta\beta}
\end{equation}
(we round off large integers).
They are defined by a Metropolis algorithm
associated with the Hamiltonian
\begin{equation}
H:\sigma\in\X\mapsto
-\frac{1}{2}\sum_{\stackrel{\scriptstyle \{i,j\}\subset\Lambda_\beta}{\scriptstyle d_1(i,j) =1}}J\sigma_i\sigma_j 
-\frac{1}{2}\sum_{i\in\Lambda} h\sigma_i
\end{equation}
for Glauber dynamics ($y\sim x$ in~(\ref{dfnmet})
has the same signification as that of the finite volume dynamics), and
\begin{equation}
H:\eta\in\X\mapsto
\sum_{\stackrel{\scriptstyle \{i,j\}\subset\Lambda_\beta}{\scriptstyle d_1(i,j) =1}}-U\eta_i\eta_j 
\end{equation}
for Kawasaki dynamics
($y\sim x$ in~(\ref{dfnmet})
now simply means that $y$ is obtained from $x=\eta\in\X$
by exchanging the value of $\eta$ between two nearest neighbour sites $i$ and $j$ in $\Lambda_\beta$).
We will assume
\begin{equation}
\Theta < \Gamma
\end{equation}
where $\Gamma$ is the energy barrier defined in the local version of the dynamics.

In large volume, we lose much of the tools
inherited from Wentzell and Freidlin.
Nevertheless many kind of restricted ensembles
have been individuated for Glauber dynamics
even in infinite volume or at fixed temperature
(\cite{DSc}, \cite{SS}).
It is not so for conservative dynamics:
as far as I know there are still only unpublished
results (\cite{GdHNOS2}, \cite{GdHNOS3})
that prove the desired properties of a set ${\cal R}$
in large volume with $d=2$.
But with the tools of potential theory
Bovier, den Hollander and Spitoni~\cite{BdHS}
computed sharp asymptotics
on some hitting times
that give very strong indication
that the set ${\cal R}$ defined as
the set of configurations for which
there are no more than $l_c(l_c-1)+1$
particles inside each square box of volume
smaller than $L^2_\beta$
with
\begin{equation}
L^2_\beta:=
e^{({\Delta}-\delta_\beta)\beta}
\mbox{ with $\delta_\beta=o(1)$ and $\frac{1}{\beta}=o(\delta_\beta)$}
\end{equation}
can be associated with a metastable restricted ensemble.
They also gave analogous results for Glauber dynamics.
I refer to \cite{BdHS} for precise statements.
Here I just want to make a few comments on the method.

The central idea is to apply~(\ref{xmeta}).
We have then four main questions to deal with:
\begin{itemize}
\item[{\bf Q1:}] How should we choose $A$ and $B$?
\item[{\bf Q2:}] How can we estimate the capacity $C_{A,B}$?
\item[{\bf Q3:}] How can we estimate the mean potential $\mu(V_{A,B})$?
\item[{\bf Q4:}] How can we link the expectation starting from the harmonic measure $E_{\nu_A}$
                 with the expectation starting from the restricted ensemble $E_{\mu_{\cal R}}$?
\end{itemize}

The ideal choices for $A$ and $B$ would be $A=\{a\}$ for any $a$ in ${\cal R}$
and $B={\cal R}^c$, then $B$ in a sequence of sets that go far away from ${\cal R}$.
If we were able to prove that for such choices $E_{\nu_A}[\tau_B]$ has a divergent
asymptotic (in exponential of $\beta$) that does not depend on $A$ and $B$,
we would not have our ``short time scale'' but, since lower bounds
on typical exit time are generally easy to get by reversibility,
the problem would essentially be solved. In particular
Q4 would have a trivial answer.
But choosing for $A$ a singleton $\{a\}$ makes in general Q2 and Q3 extremely
difficult to answer.
Indeed we always have, see~(\ref{dfncp}), 
\begin{equation}
C(a,B)\leq \mu(a)
\end{equation}
and $C(a,B)$ turns out to be super-exponentially small.
This would also imply a sharp control on $\mu(V_{A,B})$
that is extremely difficult to get:
we are indeed on the discrete version of a continuous
Dirichlet problem without solution on a bounded domain. 
In conclusion $A$ has to be ``big enough''.
As far as the choice of $B$ is concerned we remain for
a while on our ideal choice.

For a big enough $A$, Q2 is the easiest to answer:
we can make use of our two variational principles.
Actually Bovier, den Hollander and Spitoni
made use of a more elaborated variational principle
on the effective resistance that is due to
Berman and Konsowa~\cite{BK}.

Q3 is in general more difficult to answer,
since we do not have a variational principle
on the potential.
Sometimes, like in~\cite{BEGK1}, \cite{BM}, one can
reduce an estimate on a potential
to an estimate on capacities with the following
\begin{lmm}\label{pwc}
  For all $x$ outside of $A$ and $B$, disjoint sets,
  it is
  \begin{equation}
    P_x(\tau_A<\tau_B) \leq \frac{C_{x,A}}{C_{x,B}}
  \end{equation}
\end{lmm}
\noindent
{\bf Proof:} 
\begin{eqnarray}
P_x(\tau_A<\tau_B)
&=&\quad P_x(\tau_A<\tau_B|\tau_x^+>\tau_{A\cup B})P_x(\tau_x^+>\tau_{A\cup B})\nonumber\\
&&+\; P_x(\tau_A<\tau_B|\tau_x^+<\tau_{A\cup B})P_x(\tau_x^+<\tau_{A\cup B})
\end{eqnarray}
The first term in this sum is bounded from above by
\begin{equation}
P_x(\tau^+_x>\tau_A|\tau_x^+>\tau_{A\cup B})P_x(\tau_x^+>\tau_{A\cup B})
=P_x(\tau^+_x>\tau_A)
\end{equation}
while the second term is equal to
\begin{equation}
P_x(\tau_A<\tau_B)
P_x(\tau_x^+<\tau_{A\cup B})
= P_x(\tau_A<\tau_B)(1 - P_x(\tau_x^+>\tau_{A\cup B}))
\end{equation}
Solving in $P_x(\tau_A<\tau_B)$
we get
\begin{equation}
P_x(\tau_A<\tau_B)
\leq \frac{P_x(\tau^+_x>\tau_A)}{P_x(\tau_x^+>\tau_{A\cup B})}
\leq \frac{P_x(\tau^+_x>\tau_A)}{P_x(\tau_x^+>\tau_{B})}
= \frac{C_{x,A}}{C_{x,B}}
\end{equation}
\qed

\noindent
Unfortunately, single points in large volume
have not enough mass for this to be useful (see above).
Potential are difficult to estimate
(in~\cite{GdHNOS2}, \cite{GdHNOS3}
this kind of estimates involve
a complex renormalization procedure)
but we still have the trivial estimate
\begin{equation}
V_{A,B}\leq \one_{B^C}
\label{triv}
\end{equation}
On the model of our estimates of Section~\ref{sharp},
what is needed is an estimate of the kind
\begin{equation}
\mu(V_{A,B}) \sim \mu({\cal R})
\end{equation}
and this is given by~(\ref{triv})
only if $B$ is not ``too far'' from ${\cal R}$.

Denoting by $A_n$ the subset of ${\cal R}$
such that
there are no more than $n\leq l_c(l_c-1)+1$
particles inside each square box of volume
smaller than $L^2_\beta$,
by $n_0$ the smallest $n$ for which
\begin{equation}
\mu_{{\cal_R}}(A_n)\sim 1
\end{equation}
and, for all $l\geq l_c$, by $B_n$ 
the subset of ${\cal R}^c$
of the configurations $\eta$ in $\X$
that contain (as subset of $\Lambda_\beta$)
a square of side length $l$,
Bovier, den Hollander and Spitoni
give the sharp asymptotic
\begin{equation}
\forall n \in [n_0,lc(lc -1)+1],\; \forall l\in[l_c,2l_c-1],\;
E_{\nu_{A_n}}[\tau_{B_l}]\sim
\frac{3\Delta\beta e^{\Gamma\beta}}
{4\pi l_c^2(l_c^2-1)|\Lambda_\beta|} 
\label{as}
\end{equation}
Postponing the discussion on Q4,
this essentially gives us our long time scale
(beyond exponential asymptotics!)
and the constraint $l\leq 2l_c-1$
is here just as a consequence
of the difficulties that are encountered
for estimating
potentials.
However,
in the case of Glauber dynamics
and as far as I understand, 
one could remove this constraint
by an attractiveness argument
that was used in the previous works
in large volume and
makes locally available
the tools of Wentzell and Freidlin theory
(just like we used it in Section~\ref{sharp}
in alternative to Lemma~\ref{pwc}).

One of the main strength of~(\ref{as})
is that it gives asymptotics that do not depend
on $n$. One could have thought that the harmonic
measure $\nu_A$ being concentrated on the internal 
border of $A$ would have introduced a bias.
It is not so and the reason for this
is probably the same that led us from the beginning
to associate with ${\cal R}$ 
the restricted ensemble $\mu_{\cal R}$ that
is also the reversible and invariant measure
for the dynamics restricted to ${\cal R}$:
for any ``good'' ${\cal R}$ the system
should typically relax in a short time scale
to $\mu_{\cal R}$. Actually, Q4 raises a problem
of convergence to (metastable) equilibrium.
This is the object of the next and final section.

\section{Application to convergence to equilibrium}
\label{appl3}

\subsection{Spectral gap}

For an extended covering of convergence to equilibrium
we refer to~\cite{SC} and~\cite{Pe}.
Here we just indicate how the objects
we discussed above link to the argument.
We recall that for $\xi$ Markov chain
on a finite state space $\X$ with transition
probability matrix $M$,
$\xi$ is reversible with respect to the probability
measure $\mu$ if and only if $M$ is a self-adjoint
operator on $\ell^2(\mu)$.
In this case $M$ has only real eigenvalues 
\begin{equation}
1=\lambda_0 > \lambda_1 \geq \dots \geq \lambda_{N-1}\geq -1
\end{equation}
and $\lambda_{N-1} > -1$ if and only if $\xi$ is aperiodic.
In this case the rate of convergence to equilibrium
in $\ell^2(\mu)$ is governed by 
\begin{equation}
\bar\lambda:= \max\{|\lambda_1| ;|\lambda_{N-1}|\}
\end{equation} 
If $\bar\lambda\neq\lambda_1$ then the transition matrix $(M+I)/2$ 
of the associated lazy chain $\xi'$ has only positive eigenvalues
and a rate of convergence of the same order.
We will then assume that $\xi$ itself is a reversible and ergodic
Markov chain with a transition matrix $M$
the eigenvalues of which are all positive.
In this case the rate of convergence of $\xi$ in $\ell^2(\mu)$
is given by the {\em spectral gap}
\begin{equation}
\lambda:= 1 -\lambda_1
\end{equation}
that is the smallest non-zero eigenvalue
of
\begin{equation}
I-M=-\L
\end{equation}
or, equivalently,
\begin{eqnarray}
\lambda:= \min_{f\in (Ker\L)^\bot} \frac{\langle f,(-\L) f \rangle_\mu}{\langle f,f\rangle_\mu}
\label{dfnlm}
\end{eqnarray}
The numerator in this variational principle is
\begin{eqnarray}
\langle f,(-\L) f \rangle_\mu
&=&\sum_x \mu(x)f(x)\left(f(x)-\sum_y p(x,y)f(y)\right)\\
&=&\sum_x \mu(x)f(x)\sum_y p(x,y)(f(x)-f(y))\\
&=&\frac{1}{2}\sum_{x,y} \mu(x)p(x,y)(f(x)-f(y))^2\\
&=& \D(f)
\end{eqnarray}
In addition, the kernel of $\L$
is the one-dimensional subspace of $\ell^2(\mu)$
that contains all the constant functions, and the orthogonal
projection of any $f$ on this subspace is the constant function
$\mu(f)$. As a consequence one can extend the minimum
in~(\ref{dfnlm}) as a minimum on all the non-constant functions
replacing $f$ by $f-\mu(f)$, and this gives
\begin{eqnarray}
\lambda= \min_{\Var(f)\neq 0} \frac{\D(f)}{\Var(f)}
\end{eqnarray}

Any test function $f$ gives an upper bound
on the spectral gap.
Restricting the minimum to characteristic functions
we get
\begin{eqnarray}
\lambda
&\leq& \min_{A\subset\X} \frac{\sum_{e\in\partial A} c(e)}{\mu(A)(1-\mu(A))}\\
&\leq& 2\min_{\mu(A)\leq \frac{1}{2}} \frac{C(A,A^c)}{\mu(A)} \; =\; 2I
\end{eqnarray}
where the {\em isoperimetric constant} $I$ is defined by the last equation.
 
Actually $I$ gives also a {\em lower} bound on $\lambda$:

\medskip\par\noindent
{\bf Lemma [Cheeger]:}
\begin{equation}
\frac{I^2}{2}\leq \lambda \leq 2I
\end{equation}

\medskip\par\noindent
We refer to~\cite{SC} for the proof where it is shown that
a lower bound on $I$ expresses an $\ell^1$ version of a Poincar\'e inequality.
A Poincar\'e inequality is an inequality of the form
\begin{equation}
\forall f \in \ell^2(\mu),\;
\Var(f)\leq \kappa \D(f)
\end{equation}
and is equivalent to a lower bound on the spectral gap.

If instead of restricting the minimum to characteristic
function that are particular cases of equilibrium potential
we restrict the minimum
to general equilibrium potential $V_{A,B}$,
we get
\begin{eqnarray}
\lambda &\leq& \min_{A\cap B=\emptyset} \frac{C(A,B)}{\Var(V_{A,B})}\\
&\leq& \min_{A\cap B=\emptyset} \frac{C(A,B)}{\mu(A)\mu(B)}
\end{eqnarray}
Indeed,
\begin{eqnarray}
\Var(V_{A,B}) &=& \frac{1}{2}\sum_{x,y}\mu(x)\mu(y)[V_{A,B}(x)-V_{A,B}(y)]^2\\
&\geq& \frac{1}{2}\sum_{x,y\in A\cup B}\mu(x)\mu(y)\\
&=& \mu(A)\mu(B)\\
\end{eqnarray}
In the metastable situation,
when we have $a$, $b$, ${\cal G}$ that
satisfy (P1)-(P4) of Section~\ref{fv}
this gives
\begin{equation}
\lambda\leq \frac{C(A,B)}{\mu(A)\mu(B)}  
\sim \frac{C(A,B)}{\mu(A)}
\sim \frac{1}{E_{\nu_A}[\tau_B]}
\sim \frac{1}{E_a[\tau_b]}
\label{up}
\end{equation}
We prove in the next section
that this is the correct asymptotics.

\subsection{Lower bounds for the spectral gap}

We start with an easy estimate.
Given $f\in\ell^2(\mu)$ we have for any $x$ and $y$ in $\X$,
by bilinearity of $\D$ and Dirichlet's principle,
\begin{equation} 
[f(x)-f(y)]^2\leq R(x,y) D(f) 
\label{bad}
\end{equation}
Multiplying by $\mu(x)\mu(y)$ and summing on all $x$ and $y$
we get a Poincar\'e inequality
\begin{equation} 
\Var(f)\leq \frac{1}{2}\left(\sum_{x,y}\mu(x)\mu(y)R(x,y)\right) D(f) 
\end{equation}

In the metastable situation with (P1)-(P4)
one can estimate the effective resistance
between $x$ and $y$ by building a linear
(i.e., without ramifications) unitary flow from
$x$ to $y$ with $H(x)\geq H(y)$ to get
\begin{equation}
\mu(x)\mu(y)R(x,y) = O\left(\mu(y)e^{\beta\Gamma_{x,y}}\right)
\end{equation}
with $\Gamma_{x,y}$ the {\em energy barrier}
between $x$ and $y$ defined
by~(\ref{minmax}) with $x$ and $y$
in place of $a$ and $b$.
Since
\begin{equation}
\mu(a)\mu(b)R(a,b)\sim\mu(a)R(a,b)\sim E_a[\tau_b]
\end{equation}
is logarithmically equivalent to $e^\Gamma\beta$,
we get, with (P1) and (P2),
\begin{equation}
\frac{1}{\lambda}\leq E_a[\tau_b](1+o(1))
\end{equation}
and, together with~(\ref{up}),
\begin{equation}
\frac{1}{\lambda} \sim E_a[\tau_b]
\end{equation}
In \cite{BEGK1} Bovier, Eckhoff, Gayrard and Klein
prove such a relation for all the 
``low-lying eigenvalues'' of the generator.
In addition they prove that associated eigenvectors
are equivalent to some equilibrium potentials.
 
\medskip\par\noindent
{\bf Better bounds.}
In many situations our previous Poincar\'e
inequality
is however a bad one.
This is because
for different $x$ and $y$
the equality in~(\ref{bad})
is in general realised by very different $f$.
To improve this bound
we modify the network in a specific way for each $x$ and $y$.
If we increase the conductance of the edges that are
more charged by $i_{x,y}$ then we can decrease the resistance
in~(\ref{bad}). After that we will have to reconstruct a Poincar\'e
inequality on the global network on the basis of these different inequalities  
in different networks. This may give good bounds if the currents
in the different modified networks tend to use different edges.
To put it formally,
we associate with each $(x,y)$ in $\X^2$ a unitary flow $\phi_{x,y}$
and a weight function 
\begin{equation}
w_{x,y}: e\in\E \mapsto w_{x,y}(e) \in [0;+\infty[
\label{cn1}
\end{equation}
such that, for all $e$ in $\E$,
\begin{equation}
(\phi_{x,y}(e) =  0) \Rightarrow (w_{x,y}(e) = 0)
\label{cn2}
\end{equation}
Interesting choices are
\begin{eqnarray}
w^1_{x,y} &:=& 1-\one_{\{0\}}\circ \phi_{x,y}\\
w^2_{x,y} &:=& |\phi_{x,y}|\\
w^3_{x,y} &:=& |r\phi_{x,y}|\\
w^4_{x,y} &:=& r\phi^2_{x,y}
\end{eqnarray}
Then we denote by $\D_{x,y}$ the Dirichlet form
associated with the network obtained by replacing
$c$ with $cw_{x,y}$ and restriction to the connected
component that contain both $x$ and $y$.
In particular we have
\begin{equation}
\D_{x,y}(\phi_{x,y}) = \sum_{e:\phi_{x,y}(e)>0} \frac{r(e)}{w_{x,y}(e)}\phi^2_{x,y}(e)
\end{equation}
Now using both our variational principles we have
\begin{eqnarray}
|f(x)-f(y)|^2
&\leq& \D_{x,y}(\phi_{x,y}) \D_{x,y}(f)\\
& = & \D_{x,y}(\phi_{x,y}) \sum_{e:\phi_{x,y}(e)>0} c(e)w_{x,y}(e)|-\nabla_e f|^2
\end{eqnarray}
Multiplying by $\mu(x)\mu(y)$ and summing on all $x$ and $y$
we get
\begin{eqnarray}
\Var(f)
&\leq& \frac{1}{2}\sum_{x,y}\mu(x)\mu(y) \D_{x,y}(\phi_{x,y})
\sum_{\stackrel{\scriptstyle e\in\E}{\phi_{x,y}(e)>0}}\!\!\!
c(e)w_{x,y}(e)|-\nabla_e f|^2\\
&=& \frac{1}{2}\sum_{e\in\E}c(e)|-\nabla_e f|^2 
\sum_{\stackrel{\scriptstyle {x,y}\in\X}{\phi_{x,y}(e)>0}}\!\!\!
\mu(x)\mu(y) w_{x,y}(e)\D_{x,y}(\phi_{x,y})\\
&\leq& \D(f)\times
\max_{e\in\E}\!\!\!
\sum_{\stackrel{\scriptstyle x,y\in\X}{\phi_{x,y}(e)>0}}\!\!\!
\mu(x)\mu(y) w_{x,y}(e)\D_{x,y}(\phi_{x,y})
\end{eqnarray} 
This Poincar\'e inequality
gives
\begin{lmm}
  For any family of unitary flow $\phi_{x,y}$ from $x$ to $y$
  and $w_{x,y}$ that satisfies~(\ref{cn1}) and~(\ref{cn2})
  we have
  \begin{equation}
  \frac{1}{\lambda}\leq
  \max_{e\in\E}\!\!\!
  \sum_{\stackrel{\scriptstyle {x,y}\in\X}{\phi_{x,y}(e)>0}}\!\!\!
  \mu(x)\mu(y) w_{x,y}(e)\D_{x,y}(\phi_{x,y})
  \end{equation}
\end{lmm}

With $w=w^1$ we get
\begin{equation}
  \frac{1}{\lambda}\leq
  \max_{e\in\E}\!\!\!
  \sum_{\stackrel{\scriptstyle {x,y}\in\X}{\phi_{x,y}(e)>0}}\!\!\!
  \mu(x)\mu(y)\D(\phi_{x,y})
\end{equation}

With $w=w^2$ and $\phi_{x,y}$ built
as a family of linear flows that is
\begin{equation}
\phi_{x,y}(e)=\one_{\gamma(x,y)}(e)-\one_{\gamma(x,y)}(-e)
\end{equation}
with $\gamma(x,y)$ simple path from $x$ to $y$
we get
\begin{equation}
  \frac{1}{\lambda}\leq
  \max_{e\in\E}\!\!\!
  \sum_{\stackrel{\scriptstyle {x,y}\in\X}{e\in\gamma(x,y)}}
  \mu(x)\mu(y)\sum_{e'\in\gamma(x,y)} r(e')
\end{equation}
and this is Diaconis and Stroock's estimate~\cite{DSt}.

With $w=w^3$ and $\phi_{x,y}$ built like previously
we get
\begin{equation}
  \frac{1}{\lambda}\leq
  \max_{e\in\E}\!\!\!
  \sum_{\stackrel{\scriptstyle x,y\in\X}{e\in\gamma(x,y)}}
  \mu(x)\mu(y)r(e)|\gamma(x,y)|
\end{equation}
that is one of Sinclair's estimates in~\cite{Si},
while in the general case
we get 
\begin{equation}
  \frac{1}{\lambda}\leq
  \max_{e\in\E}\!\!\!
  \sum_{\stackrel{\scriptstyle {x,y}\in\X}{\phi_{x,y}(e)>0}}\!\!\!
  \mu(x)\mu(y)r(e)\phi_{x,y}(e)
  \sum_{\stackrel{\scriptstyle e'\in\E}{\phi_{x,y}(e')>0}}\!\!\!
  \phi_{x,y}(e')
\end{equation}
that is essentially equivalent to Sinclair's formula
for ``multicommodity flows'', actually strictly equivalent
if $\phi_{x,y}$ is built (like usually it is)
as a flow of geodesic paths.

With $w=w^4$ 
we get
\begin{equation}
  \frac{1}{\lambda}\leq
  \max_{e\in\E}\!\!\!
  \sum_{\stackrel{\scriptstyle x,y\in\X}{\phi_{x,y}(e)>0}}\!\!\!
  \mu(x)\mu(y)r(e)\phi^2_{x,y}(e)
  \Big|\Big\{e'\in\E :\: \phi_{x,y}(e')>0\Big\}\Big|
\end{equation}
Previous choices have proven their efficiency
but at the time of sending these notes
I did not have that of seriously checking
the advantages of this last choice.

\subsection{Mixing time}

A stronger notion of convergence to equilibrium
is that of the {\em mixing time}
\begin{equation}
\tau_1
:= \inf \left\{t\geq0 :\: \max_{x\in\X}
\left\|P_x^t - \mu\right\|_{TV} 
\leq \frac{1}{e}\right\}
\end{equation}
where
\begin{equation}
\left\|P_x^t - \mu\right\|_{TV} 
:= \max_{A\in\X}\Big|P_x(\xi(t)\in A)-\mu(A)\Big|
\end{equation}
is the total variation distance
that, maximized over $x$, decays exponentially~\cite{Pe}.
Since an exponential decay of the total variation
distance implies an exponential decay at the same rate (at least)
in $\ell^2(\mu)$, we have
\begin{equation}
-\ln(1-\lambda) \geq \frac{1}{\tau_1}
\end{equation}

In our metastable situation,
since the equilibrium measure
is essentially concentrated on $\{b\}$,
the study of the convergence to equilibrium
should be essentially an hitting time problem.
The coupling argument used in~\cite{OV},
turns this intuition in a rigorous estimate
of the mixing time.
For any ergodic and aperiodic Markov chain
$\xi$ and any coupling $(X,Y)$
with marginals following
the law of $\xi$ we have, for all $t\geq 0$~\cite{Pe}
\begin{equation}
\max_{x\in\X}
\left\|P_x^t - \mu\right\|_{TV} 
\leq 
\max_{x,y\in\X}
P_{x,y}(\tau_c>t)
\label{cmix}
\end{equation}
with
\begin{equation}
\tau_c:=\inf\{t\geq 0 :\: X(t)=Y(t)\}
\end{equation}
In our metastable situation,
we can then obtain sharp estimates
on the basis of the {\em exponential law}:
\begin{prp}
  For a continuous time Markov chain
  built on a Metropolis algorithm
  at low temperature $\beta^{-1}$
  on a finite state space $\X$ and such that
  (P1)-(P4) hold, 
  the random variables
  \begin{equation}
  \theta_b:=\frac{\tau_b}{E[\tau_b]}
  \quad\mbox{and}\quad
  \theta_B:=\frac{\tau_B}{E[\tau_B]}
  \end{equation}  
  converge in law under $P_a$ 
  towards an exponential random time
  of mean 1 as $\beta$ goes to infinity.
\end{prp}

\noindent
This kind of result goes back at least
to Cassandro, Galves, Olivieri and Vares' work
\cite{CGOV}.
The nice following proof
is adapted from~\cite{MNOS}.
The argument works the same
in much more general situation:
as long as one can talk about metastable single {\em states}.

\noindent
{\bf Proof:}
We prove the result for $\theta_B$:
the same proof works
for $\theta_b$.
Since 
\begin{equation}
T\geq 0\mapsto P_a(\tau_B > T )
\end{equation}
is a non-increasing continuous function
taht goes from 1 to 0,
we can define
\begin{equation}
T_B:=\min\left\{T\geq 0:\: P_a(\tau_B > T ) = \frac{1}{e}\right\}
\label{dfnt}
\end{equation}
For all positive $s$ we define also
\begin{equation}
\tau^*:=\inf\Big\{t\geq sT_B :\: X(t)\in \{a\}\cup B\Big\}
\end{equation}
With
\begin{equation}
R:= e^{(\tilde\Gamma+\delta_0)\beta}
\end{equation}
where $\tilde\Gamma$ is the internal resistance
of $B^c$, that is the maximal energy barrier
under $\Gamma$ between two states $x$ in $B^c$
and $y$ in $B^c \cup \partial_+(B^c)$, and $\delta_0>0$ is such that $\tilde\Gamma+\delta_0<\Gamma$,
we know, by classical Wentzell-Freidlin
theory, that
the probability of $\tau^*-sT_B$ being
larger than $R$ is super exponentially small
($SES$), that is
\begin{equation}
\limsup_{\beta\rightarrow +\infty}
\frac{1}{\beta}\ln P_a(\tau^*-sT_B \geq R)
= -\infty
\end{equation}
For all positive $t$, if $\tau_B$ is larger than $(s+t)T_B$ then
$X$ will be outside $B$ on the whole intervals
\begin{equation}
I=[0, sT_B] 
\quad\mbox{and}\quad 
J=[\tau^*\wedge (s+t)T_B, (s+t)T_B]
\end{equation}
Since, up to a $SES$ event, $J$ is longer than $tT_B - R$,
\begin{equation}
P_a(\tau_B>(s+t)T_B)\leq 
P_a(\tau_B>sT_B)
P_a(\tau_B>tT_B-R)
+SES
\label{maj}
\end{equation}
If $X$ stays outside  
$B$ on  
\begin{equation}
I=[0, \tau^*] \quad\mbox{and}\quad  J=[\tau^*,\tau^*+tT_B]
\end{equation}
then $\tau_B$ will be larger than $(s+t)T_B$.
Since, up to a $SES$ event, $I$ is shorter than $tT_B + R$,
\begin{equation}
P_a(\tau_B>sT_B+R)
P_a(\tau_B>tT_B)
-SES
\leq 
P_a(\tau_B>(s+t)T_B)
\label{min}
\end{equation}
By classical Wentzell-Freidlin theory
we have also, for all $\delta>0$,
\begin{equation}
\limsup_{\beta\rightarrow+\infty}
\frac{1}{\beta}\ln
P_a (\tau_B < e^{(\Gamma-\delta)\beta})
< 0
\end{equation}
so that
\begin{equation}
R=o(T_B)
\label{o}
\end{equation}
and, by~(\ref{maj}), for all large enough $\beta$,
\begin{equation}
P_a(\tau_B>(s+t)T_B)\leq 
P_a(\tau_B>sT_B)
P_a(\tau_B>(t-1)T_B)
+SES
\end{equation}
For all $n=2k$, we get then, by induction
on $k\geq 2$ and using~(\ref{dfnt})
\begin{equation}
P_a(\tau_B>nT_B)\leq 
e^{-\frac{n}{2}\beta} + n SES 	
\end{equation}
and the tightness of $(\bar\theta_\beta)_{\beta>0}$
with
\begin{equation}
\bar\theta_\beta:=\frac{\tau_B}{T_B}
\end{equation}
By Alaoglu's theorem and diagonal extraction
there is a diverging sequence $(\beta_k)_{k\geq 0}$
such that $\bar\theta_{\beta_k}$ converges in law
towards a random variable $\hat\theta$,
for which, with~(\ref{maj}), (\ref{min}) and (\ref{o}),
we have necessarily 
\begin{equation}
P\left(\hat\theta > s+t\right) = P\left(\hat\theta>s\right)
P\left(\hat\theta>t\right)
\end{equation}
for all $s$ and $t$ outside ${\cal A}$, set of $\hat\theta$'s
atoms.
Since ${\cal A}$ is a countable set, so is $\Q{\cal A}$
and there is some $x>0$ such that
\begin{equation}
x\Q\cap{\cal A}\subset \{0\}
\end{equation}
By density of $x\Q$ and monotonicity of the repartition
function
we conclude that $\hat\theta$ follows an exponential law.
By~(\ref{dfnt}) this can only be that of mean 1
and we get, for all $t>0$,
\begin{equation}
  \lim_{\beta\rightarrow +\infty}
  P_a\Big(\tau_B > tT_B\Big)
  =e^{-t}
\end{equation}
In addition, by dominated convergence,
\begin{equation}
\lim_{\beta\rightarrow+\infty}\frac{E_a[\tau_B]}{T_B}
= \lim_{\beta\rightarrow+\infty}\int_0^{+\infty}P_a(\bar\theta_\beta > t) dt
= \int_0^{+\infty}e^{-t} dt
= 1
\end{equation}
Now for $t>0$ and any $\epsilon>0$ we choose $\gamma>1$
such that
\begin{equation}
e^{-t}-\epsilon 
< \frac{1}{\gamma}e^{-\gamma t}
< {\gamma}e^{-\frac{t}{\gamma}} 
< e^{-t}+\epsilon 
\end{equation}
For $\beta$ large enough
we have, on the one hand 
\begin{equation}
P_a(\theta_B>t)
\leq P_a \left(\tau_B > \frac{tT_B}{\gamma}\right)
\leq {\gamma}e^{-\frac{t}{\gamma}}
\leq e^{-t}+\epsilon 
\end{equation}
on the other hand
\begin{equation}
P_a(\theta_B>t)
\geq P_a \left(\tau_B > {\gamma}{tT_B}\right)
\geq \frac{1}{\gamma}e^{-\gamma {t}} 
\geq e^{-t}-\epsilon 
\end{equation}
and this concludes the proof.
\qed

\medskip\par
A simple analysis of the basic coupling, that is $X$ coupled with an independent Markov chain $Y$
with same generator, is then sufficient to give a sharp estimate of the mixing time.
If $X$ starts from $x$ and $Y$ starts from $y$, then
in a time shorter than $e^{(\Gamma-\delta)\beta}$, for a small enough $\delta>0$,
both will have typically reached $A\cup B$.
If they have reached the same cycle they will also have met with a probability
exponentially close to one (look at the dynamics
as one built on a Metropolis algorithm on the product space $\X\times\X$).
If they have reached different cycles,
say $X$ reached $A$ and $Y$ reached $B$,
then $Y$ will typically stay in $B$
for a time larger than $e^{(\Gamma+\delta)\beta}$
(for $\delta$ small enough),
and if $X$ reached $B$ before $Y$ left it
the two processes will typically meet
in a time shorter than $e^{(\Gamma -\delta)\beta}$.
The key question is then that of the distribution
of the hitting time of $B$ for the process $X$
that reached $A$ (and $a$) on the short time scale 
$e^{(\Gamma -\delta)\beta}$.
The exponential law provides the answer.
Putting everything together we get that
for any $\gamma >1$ there is a positive $\beta_0$
such that, for all $\beta$ larger than $\beta_0$,
\begin{equation}
\max_{x,y}P_{x,y}\left(\tau_c>\gamma E_a[\tau_B]\right)
\leq
\frac{1}{e}
\end{equation}
With~(\ref{cmix}) this gives
\begin{equation}
\tau_1\leq \gamma E_a[\tau_B]
\end{equation}
or, for any  $\gamma>1$ and $\beta$ large enough,
\begin{equation}
\tau_1\leq \gamma E_a[\tau_b]
\end{equation}
Since we also have
\begin{equation}
E_a[\tau_b]
\sim \frac{1}{\lambda}
\sim \frac{-1}{\ln(1-\lambda)}
\leq \tau_1
\end{equation}
we conclude
\begin{equation}
\tau_1
\sim 
E_a[\tau_b]
\end{equation}

\subsection{A final remark}

All this gives the tools to deal with
relaxation to the {\em global} equilibrium
$\mu$.
But this is not conclusive
as far as the question raised
at the end of Section~\ref{appl2}
is concerned: what about the relaxation
to the restricted ensemble, or to a {\em local}
equilibrium?
As far as I understand what has been done
is the last decades in the study of metastability
it seems to me that this is now a central
and still unsolved question.
One source of inspiration to deal with it
could be found in Miclo's work~\cite{MI}.

\end{document}